\documentclass[12pt]{article}
\usepackage{amsmath, amsthm,
amssymb,extsizes}
\usepackage[all]{xy}
\DeclareMathOperator{\Aut}{Aut}
\DeclareMathOperator{\End}{End}
\DeclareMathOperator{\ind}{ind}
\DeclareMathOperator{\rep}{Rep}
\DeclareMathOperator{\Isom}{Is}
\DeclareMathOperator{\charact}{char}
\DeclareMathOperator{\diag}{diag}

\title{Canonical
matrices of isometric
operators on
indefinite inner
product spaces\footnotetext{This is the author's version of a work that was accepted for publication in Linear Algebra and its Applications (2007), doi:10.1016/j.laa.2007.08.016.}}
\author{
Vladimir V. Sergeichuk%
\thanks{The research was done
while the author was
visiting the
University of S\~ao
Paulo supported by
FAPESP, processo
05/59407-6.}
\\
Institute of
Mathematics,
Tereshchenkivska 3,
Kiev,
Ukraine,\\sergeich@imath.kiev.ua}
\date{}

\begin{document}
 \maketitle

\newcommand{\ddd}{
\text{\begin{picture}(12,8)
\put(-2,-4){$\cdot$}
\put(3,0){$\cdot$}
\put(8,4){$\cdot$}
\end{picture}}}

\newcommand{\dia}{\,\diagdown\,}

\newcommand{\lin}{\;\frac{}{
\phantom{\longrightarrow}}\;}

\newcommand{\ci}{
\begin{picture}(6,6)
\put(3,3){\circle*{3}}
\end{picture}}

\newcommand{\is}{\stackrel
{\text{\raisebox{-1ex}{$\sim\
\;$}}}{\to}}

\newtheorem{theorem}{Theorem}[section]
\newtheorem{lemma}{Lemma}[section]
\newtheorem*{corollary}{Corollary}

\theoremstyle{remark}
\newtheorem*{remark}{Remark}
\newtheorem*{example}{Example}

\renewcommand{\le}{\leqslant}
\renewcommand{\ge}{\geqslant}

\begin{abstract}
We give canonical
matrices of a pair
$(A,B)$ consisting of
a nondegenerate form
$B$ and a linear
operator $A$
satisfying
$B(Ax,Ay)=B(x,y)$ on a
vector space over
$\mathbb F$ in the
following cases:
\begin{itemize}
  \item
$\mathbb F$ is an
algebraically closed
field of
characteristic
different from $2$ or
a real closed field,
and $B$ is symmetric
or skew-symmetric;

  \item
$\mathbb F$ is an
algebraically closed
field of
characteristic $0$ or
the skew field of
quaternions over a
real closed field, and
$B$ is Hermitian or
skew-Hermitian with
respect to any
nonidentity involution
on $\mathbb F$.
\end{itemize}
These classification
problems are wild if
$B$ may be degenerate.

We use a method that
admits to reduce the
problem of classifying
an arbitrary system of
forms and linear
mappings to the
problem of classifying
representations of
some quiver. This
method was described
in [V.V. Sergeichuk,
{\it Math. USSR-Izv.}
31 (1988) 481--501].


 \end{abstract}

\section{Introduction}
\label{s_intr}

Let $\mathbb F$ be a
field or skew field of
characteristic
different from 2 with
involution (which may
be the identity). We
consider the problem
of classifying pairs
\begin{equation}\label{aa}
(A,B)
\end{equation}
consisting of a
nondegenerate
Hermitian or
skew-Hermitian form
$B\colon V\times
V\to\mathbb F$ on a
right vector space $V$
over $\mathbb F$ and
an operator $A\colon
V\to V$ that is
\emph{isometric} with
respect to $B$; i.e.,
\begin{equation*}\label{ybt}
 B(Au,Av)=B(u,v)\qquad\text{for all }
u,v\in V.
\end{equation*}

This problem was
solved in
\cite[Theorem
5]{ser_izv} over
$\mathbb F$ up to
classification of
Hermitian forms over
finite extensions of
$\mathbb F$, we
present the solution
in Theorem
\ref{Theorem 5}.  This
implies its complete
solution over $\mathbb
C$ and $\mathbb R$
since the
classification of
Hermitian forms over
$\mathbb C$ and
$\mathbb R$ is known.
But the canonical
matrices in
\cite{ser_izv} are not
simple since they are
based on the Frobenius
canonical form over
$\mathbb F$ for
similarity.

\emph{The first
purpose} of this paper
is to give simple
canonical matrices of
pairs \eqref{aa} over
an algebraically or
real closed field of
characteristic
different from $2$
basing on the Jordan
canonical form for
similarity. We also
obtain canonical
matrices of $(A,B)$
over the skew field
$\mathbb H$ of real
quaternions; they are
given in \cite{ser1}
incorrectly (see the
footnote on page
\pageref{page}). This
classification problem
was studied in
\cite{huppert1,
huppert2,miln,sch},
other canonical
matrices of $(A,B)$
and their applications
are given in
\cite{au-Rod,rod_pert}
over $\mathbb C$ and
$\mathbb R$, and in
\cite{apl-Rod} over
$\mathbb H$.

\emph{The second
purpose} of this paper
is to present in
sufficient detail a
technique for
classifying systems of
forms and linear
mappings (we use it to
obtain canonical
matrices of
\eqref{aa}). It was
devised by Roiter
\cite{roi} and the
author
\cite{ser_first,
ser_disch,ser_izv}. It
is practically unknown
although many
classification
problems solved
recently could be
easily solved by this
method. This
linearization
technique reduces the
``nonlinear'' problem
of classifying an
arbitrary system $\cal
S$ of forms and linear
mappings over a field
or skew field $\mathbb
F$ of characteristic
different from $2$
\begin{itemize}
  \item
to the ``linear''
problem of classifying
some system
$\underline{\cal S}$
of linear mappings
over $\mathbb
F$---i.e., to the
problem of classifying
representations of a
quiver with relations,
and
  \item
to the problem of
classifying Hermitian
forms over fields or
skew fields that are
finite extensions of
the center of $\mathbb
F$.
\end{itemize}

The corresponding
reduction theorems
were extended in
\cite{ser_izv} to the
problem of classifying
selfadjoint
representations of a
linear category with
involution and in
\cite{ser_sym} to the
problem of classifying
symmetric
representations of an
algebra with
involution. Similar
theorems were proved
for bilinear and
sesquilinear forms by
Gabriel, Riehm, and
Shrader-Frechette
\cite{gab2,riehm1,riehm2};
for additive
categories with
quadratic or Hermitian
forms on objects by
Quebbermann, Scharlau,
and Schulte
\cite{que-scha,w.schar};
for generalizations of
quivers involving
linear groups by
Derksen, Shmelkin, and
Weyman
\cite{der-wey,shme}.

Two cases are possible
for the system $S$.
\medskip

\emph{Case 1:
$\underline{\cal S}$
is wild.} This means
that the problem of
classifying the system
$\underline{\cal S}$
contains the problem
of classifying pairs
of matrices up to
simultaneous
similarity. The latter
problem is hopeless
since it contains the
problem of classifying
an \emph{arbitrary}
system of linear
mappings
\cite[Theorems 4.5 and
2.1]{bel-ser_complexity}.
Hence, the problem of
classifying the system
${\cal S}$ is hopeless
too. For example, the
wildness of
$\underline{\cal S}$
was proved in
\cite[Theorems 5.4 and
5.5]{ser_prep} for the
problems of
classifying
\begin{itemize}

 \item[$-$]
selfadjoint/metric
operators on a space
with degenerate
indefinite scalar
product (we replicate
this result in Theorem
\ref{th_dic}; this
classification problem
was considered in
\cite{mehl-rod_deg})
and

  \item[$-$]
normal operators on a
space with degenerate
indefinite scalar
product (this problem
was posed in \cite[p.
84]{g-l-r}; its
wildness was also
proved in \cite{goh}).
\end{itemize}
Thus, these problems
are hopeless, and so
the problem of
classifying \eqref{aa}
cannot be solved if
$B$ may be degenerate.
\medskip

\emph{Case 2:
$\underline{\cal S}$
is not wild.} Then the
problem of classifying
the system
$\underline{\cal S}$
can be solved. In each
dimension, the set of
Belitski\u \i's
canonical matrices of
the system
$\underline{\cal S}$
consists of a finite
number of matrices and
1-parameter families
of matrices and is
presented by a finite
number of points and
straight lines in the
affine matrix space
(see \cite[Theorem
3.1]{ser-can} and also
\cite{gab_vos}). For
example, the system
$\underline{\cal S}$
is not wild for the
problems of
classifying
\begin{itemize}
  \item
sesquilinear forms,

  \item
pairs of forms, in
which the first form
is
$\varepsilon$-Hermitian
and the second is
$\delta$-Hermitian
($\varepsilon,\delta
\in\{1,-1\}$), and
  \item
isometric or
selfadjoint operators
on a space with
nondegenerate
$\varepsilon$-Hermitian
form (an operator $A$
is \emph{selfadjoint}
with respect to a form
$B$ if
$B(Ax,y)=B(x,Ay)$).
\end{itemize}
Their canonical
matrices were obtained
by the linearization
technique in
\cite{ser_disch,ser_prep}
and also in
\cite[Theorems
3--6]{ser_izv} over
any field of
characteristic
different from 2 up to
classification of
Hermitian forms over
its finite extensions.
\medskip

Theorem \ref{tetete}
implies that each
system of forms and
linear mappings over
$\mathbb C$, $\mathbb
R$, or $\mathbb H$
decomposes into a
direct sum of
indecomposable systems
uniquely up to
isomorphism of
summands. Hence, it
suffices to classify
only indecomposable
systems.

A detailed exposition
of the theory of
operators on spaces
with indefinite scalar
product is given in
the books
\cite{g-l-r,g-l-r1}.

The paper is organized
as follows. In Section
\ref{s_intre} we
formulate Theorem
\ref{theor} about
canonical matrices of
pairs \eqref{aa} over
algebraically or real
closed fields and skew
fields of real
quaternions. We also
formulate Theorem
\ref{Theorem 5}, which
is a useful
generalization of
{\cite[Theorem
5]{ser_izv}} and gives
canonical matrices of
\eqref{aa} over any
field of
characteristic
different from 2 up to
classification of
Hermitian forms.

Section \ref{s_pos}
contains a detailed
description of the
linearization
technique; it can be
read independently of
Section \ref{s_intre}.
Theorem \ref{tetete}
in this section
extents Sylvester's
Inertia Theorem to
systems of forms and
linear mappings.

In Sections
\ref{s_pro} and
\ref{secmat} we prove
Theorems \ref{theor}
and \ref{Theorem 5}.

In Section
\ref{metric} we
present Theorem 5.4 of
\cite{ser_prep} about
the wildness of the
problem of classifying
pairs \eqref{aa} in
which $B$ may be
degenerate.

\section{Canonical
matrices of isometric
operators}
\label{s_intre}

We recall some
properties of
algebraically or real
closed fields and skew
fields of real
quaternions, and
formulate Theorems
\ref{theor} and
\ref{Theorem 5} about
canonical matrices of
pairs \eqref{aa}.

\subsection{Isometric
operators over an
algebraically or real
closed field and over
quaternions}
\label{sub_01}

In this paper,
$\mathbb F$ denotes a
field or skew field of
characteristic
different from 2 with
involution $a\mapsto
\bar{a}$; that is, a
bijection $\mathbb
F\to \mathbb F$
satisfying
\[
 \overline{a+b}=\bar a+ \bar
b,\qquad
 \overline{ab}=\bar b \bar
a,\qquad
 \bar{\bar a}=a.
\]
Therefore, the
involution can be the
identity only if
$\mathbb F$ is a
field. All vector
spaces are assumed to
be finite dimensional
right vector spaces.

A mapping $B\colon
U\times V\to \mathbb
F$ on vector spaces
$U$ and $V$ over
$\mathbb F$ is called
a {\it sesquilinear
form} if
\begin{gather*}
  B(ua+u'a',v)=
  \bar{a}B(u,v)+
  \bar{a'}B(u',v),\\
  B(u,va+v'a')
  =B(u,v)a+B(u,v')a'
\end{gather*}
for all $u,u'\in U$,
$v,v'\in V$, and
$a,a'\in \mathbb F$.
This form is
\emph{bilinear} if
$\mathbb F$ is a field
and the involution
$a\mapsto\bar{a}$ is
the identity (we
consider bilinear
forms as a special
case of sesquilinear
forms). If
$e_1,\dots,e_m$ and
$f_1,\dots, f_n$ are
bases of $U$ and $V$,
then
\begin{equation*}\label{1.3}
  B(u,v)=[u]_e^{*}B_{ef}[v]_f
  \qquad\text{for all
$u\in U$ and $v\in
V$},
\end{equation*}
where $[u]_e$ and
$[v]_f$ are the
coordinate vectors,
$[u]_e^{*}:=\overline{[u]}_e^{T}$,
and $
B_{ef}:=[B(e_i,f_j)]$
is the matrix of $B$.

Let $\varepsilon$ be
an element of the
center $C(\mathbb F)$
of $\mathbb F$ such
that
$\varepsilon\bar\varepsilon=1$.
A sesquilinear form
$B\colon V\times V\to
\mathbb F$ is called
$\varepsilon$-{\it
Hermitian} if
\begin{equation*}\label{1.4}
B(u,v)=\varepsilon
\overline{B(v,u)}\qquad\text{for
all } u,v\in V;
\end{equation*}
it is called {\it
Hermitian} if
$\varepsilon=1$ and
{\it skew-Hermitian}
if $\varepsilon=-1$.
Clearly,
$\varepsilon=\pm 1$ if
the involution acts
identically on
$C(\mathbb F)$.
Without loss of
generality, \emph{we
will assume that
$\varepsilon=1$ if the
involution acts
nonidentically on
$C(\mathbb F)$} since
then an
$\varepsilon$-{Hermitian}
form $B$ can be made
Hermitian by
multiplying it by
$1+\bar\varepsilon$ if
$\varepsilon\ne -1$
because
\[
(1+\bar\varepsilon)B(u,v)=
(1+\bar\varepsilon)
\varepsilon
\overline{B(v,u)}=
(1+\varepsilon)
\overline{B(v,u)}=
\overline{(1+
\bar\varepsilon)
B(v,u)},
\]
and by $a-\bar a$ for
any $a\ne\bar a$ from
$C(\mathbb F)$ if
$\varepsilon=-1$.

Let $(A,B)$ be a pair
consisting of a
nondegenerate
$\varepsilon$-Hermitian
form $B$ and an
isometric operator $A$
on a vector space $V$.
Their matrices $A_e$
and $B_e$ in a basis
of $V$ satisfy the
conditions:
\begin{equation}\label{ugr}
B_e=\varepsilon
B_e^*=A_e^*B_eA_e,\qquad
\text{$A_e$ and $B_e$
are nonsingular,}
\end{equation}
where
$A_e^*:=\bar{A_e}^T$
(usually the letter
$H$ is used instead of
$B_e$, then $A_e$
satisfying \eqref{ugr}
is called
\emph{$H$-unitary},
see \cite{au-Rod}).
Every change of the
basis reduces
$(A_e,B_e)$ by
transformations
\begin{equation}\label{aaa}
 (A_e,B_e)\mapsto
 (S^{-1}A_eS,S^*B_eS),
 \qquad S\text{ is
nonsingular}.
\end{equation}
In Theorem \ref{theor}
we give canonical
matrices of pairs
$(A_e,B_e)$ satisfying
\eqref{ugr} with
respect to
transformations
\eqref{aaa} over:
\begin{itemize}
  \item
an {\it algebraically
closed field of
characteristic
different from} 2,
  \item
a {\it real closed
field}---i.e, a field
whose algebraic
closure has a finite
degree $\ne 1$ (see
Lemma \ref{l00}), and
  \item
the {\it skew field of
quaternions}
\begin{equation*}\label{1ya}
 {\mathbb
 H}=\{a+bi+cj+dk\,|\,a,b,c,d\in\mathbb
 P\}
\end{equation*}
over a real closed
field $\mathbb P$,
where
$i^2=j^2=k^2=-1$,
$ij=k=-ji$,
$jk=i=-kj$, and
$ki=j=-ik.$
\end{itemize}
Without loss of
generality we can
consider only two
involutions on
$\mathbb H$:
\emph{quaternionic
conjugation}
\begin{equation}\label{ne}
 a+bi+cj+dk\
\longmapsto\
 a-bi-cj-dk
\end{equation}
and \emph{quaternionic
semiconjugation}
\begin{equation}\label{nen}
 a+bi+cj+dk\
\longmapsto\
a-bi+cj+dk,
 \qquad a,b,c,d\in\mathbb P,
\end{equation}
because by Lemma
\ref{l1a} if an
involution on $\mathbb
H$ is not quaternionic
conjugation then it
has the form
\eqref{nen} in a
suitable set of
imaginary units
$i,j,k$.

There is a natural
one-to-one
correspondence
\begin{equation*}\label{1ye}
\left\{\parbox{5cm}{$\
$algebraically closed fields\\
with nonidentity
involution}\right\}
\quad\longleftrightarrow
\quad
\bigl\{\text{real
closed fields}\bigr\}
\end{equation*}
sending an
algebraically closed
field with nonidentity
involution to its
fixed field. This
follows from our next
lemma, in which we
collect known results
about such fields.

\begin{lemma}\label{l00}
{\rm(a)} Let\/
$\mathbb P$ be a real
closed field and\/ let
$\mathbb K$ be its
algebraic closure.
Then $\charact{\mathbb
P}=0$ and
\begin{equation}\label{1pp}
\mathbb K={\mathbb
P}+{\mathbb P}i,\qquad
i^2=-1.
\end{equation}
The field\/ ${\mathbb
P}$ has a unique
linear ordering $\le$
such that
\begin{equation*}\label{slr}
\text{$a>0$ and\,
$b>0$}
 \quad\Longrightarrow\quad
\text{$a+b>0$ and\,
$ab>0$}.
\end{equation*}
The positive elements
of\/ $\mathbb P$ with
respect to this
ordering are the
squares of nonzero
elements.

{\rm(b)} Let\/
$\mathbb K$ be an
algebraically closed
field with nonidentity
involution. Then
$\charact\mathbb K=0$,
\begin{equation}\label{123}
\mathbb
P:=\bigl\{k\in{\mathbb
K}\,\bigr|\,
\bar{k}=k\bigr\}
\end{equation}
is a real closed
field,
\begin{equation}\label{1pp11}
\mathbb K={\mathbb
P}+{\mathbb P}i,\qquad
i^2=-1,
\end{equation}
and the involution has
the form
\begin{equation}\label{1ii}
\overline{a+bi}=a-bi,\qquad
a,b\in\mathbb P.
\end{equation}

{\rm(c)} Every
algebraically closed
field $\mathbb F$ of
characteristic $0$
contains at least one
real closed subfield.
Hence, $\mathbb F$ can
be represented in the
form \eqref{1pp11} and
possesses the
involution
\eqref{1ii}.
\end{lemma}

\begin{proof} (a)
Let $\mathbb K$ be the
algebraic closure of
$\mathbb F$ and
suppose
$1<\dim_{\mathbb
P}{\mathbb K}<\infty$.
By Corollary 2 in
\cite[Chapter VIII, \S
9]{len}, we have
$\charact{\mathbb
P}=0$ and \eqref{1pp}.
The other statements
of part (a) follow
from Proposition 3 and
Theorem 1 in
\cite[Chapter XI, \S
2]{len}.

(b) If $\mathbb K$ is
an algebraically
closed field with
nonidentity involution
$a\mapsto \bar{a}$,
then this involution
is an automorphism of
order 2. Hence
${\mathbb K}$ has
degree $2$ over its
\emph{fixed field}
${\mathbb P}$ defined
in \eqref{123}. Thus,
${\mathbb P}$ is a
real closed field. Let
$i\in \mathbb K$ be
such that $i^2=-1$. By
(a), every element of
${\mathbb K}$ is
uniquely represented
in the form $k=a+bi$
with $a,b\in{\mathbb
P}$. The involution is
an automorphism of
${\mathbb K}$, so
$\bar{i}^2=-1$. Thus,
$\bar{i}=-i$ and the
involution has the
form \eqref{1ii}.

(c) This statement is
proved in \cite[\S 82,
Theorem 7c]{wan}.
\end{proof}

For each real closed
field, we denote by
$\le$ the ordering
from Lemma
\ref{l00}(a). Let
$\mathbb K=\mathbb
P+\mathbb Pi$ be an
algebraically closed
field with nonidentity
involution represented
in the form
\eqref{1pp11}. By the
\emph{absolute value}
of $k=a+bi\in\mathbb
K$ ($a,b\in\mathbb P)$
we mean a unique
nonnegative real root
of $a^2+b^2$, which we
write as
\begin{equation}\label{1kk}
|k|:=\sqrt{a^2+b^2}
\end{equation}
(this definition is
unambiguous since
$\mathbb K$ is
represented in the
form \eqref{1pp11}
uniquely up to
replacement of $i$ by
$-i$). For each
$M\in{\mathbb
K}^{m\times n}$, its
{\it realification}
$M^{\mathbb
P}\in{\mathbb
P}^{2m\times 2n}$ is
obtained by replacing
every entry $a+bi$ of
$M$ by the $2\times 2$
block
\begin{equation}\label{1j}
\begin{matrix}
 a&-b\\b&a
\end{matrix}
\end{equation}

Define the $n\times n$
matrices
\begin{equation}\label{ase}
J_n(\lambda) :=
\begin{bmatrix}
  \lambda&1&&0\\
  &\lambda&\ddots&\\
  &&\ddots&1\\
0&&&\lambda
\end{bmatrix},\qquad
\Lambda_n:=\begin{bmatrix}
1&2&\cdots&2
\\&1&\ddots&\vdots
\\&&\ddots&2
\\0&&&1
\end{bmatrix},
\end{equation}
\begin{equation}\label{ases}
F_n:=\begin{bmatrix}
0&&&&\ddd&
\\&&&-1
\\&&1&
\\&-1&\\
1&&&&0
\end{bmatrix}.
\end{equation}

If $M$ is nonsingular
it is convenient to
write
\begin{equation*}\label{l.dedk}
M^{-*}:=(M^{-1})^*,\qquad
M^{-T}:=(M^{-1})^T.
\end{equation*}
The {\it skew sum} of
two matrices is
defined by
\begin{equation*}\label{1.2a}
M\dia N:=
\begin{bmatrix}0&N\\M &0
\end{bmatrix}.
\end{equation*}

The main result of
this paper is the
following theorem.

\begin{theorem}
\label{theor}

Let $\mathbb F$ be one
of the following
fields or skew fields:
\begin{itemize}
  \item[\rm(a)]
an algebraically
closed field of
characteristic
different from $2$
with the identity
involution;

  \item[\rm(b)]
an algebraically
closed field with
nonidentity
involution;

  \item[\rm(c)]
a real closed field
$\mathbb P$ $($by
Lemma \ref{l00}, its
algebraic closure is
represented in the
form\/ $\mathbb
P+\mathbb P i$ and
possesses the
involution
$a+bi\mapsto a-bi)$;

  \item[\rm(d)]
the skew field\/
$\mathbb H=\mathbb
P+\mathbb P i+ \mathbb
P j+\mathbb P k$ of
quaternions over a
real closed field\/
$\mathbb P$, with
quaternionic
conjugation \eqref{ne}
or quaternionic
semiconjugation
\eqref{nen}.
\end{itemize}
Let $\varepsilon =\pm
1$ $(\varepsilon =1$
if $\mathbb F$ is
{\rm(b))} and let
$(A,B)$ be a pair
consisting of a
nondegenerate
$\varepsilon
$-Hermitian form $B$
on a right vector
space over $\mathbb F$
and an operator $A$ on
this space that is
isometric with respect
to $B$.

Then there exists a
basis in which $(A,B)$
is given by a direct
sum, determined
uniquely up to
permutation of
summands,
respectively,

\begin{itemize}
\item[\rm(a)] of the
following matrix pairs
that are given by
$0\ne\lambda\in\mathbb
F$ determined up to
replacement by
$\lambda^{-1}$:

\begin{itemize}
  \item[\rm(i)]
$(J_n(\lambda)\oplus
J_n(\lambda)^{-T},
I_n\dia \varepsilon
I_n)$, except  for
$\lambda=\pm 1$ and
$\varepsilon
=(-1)^{n+1}$,

  \item[\rm(ii)]
$(\lambda\Lambda_n,F_n)$
if $\lambda=\pm 1$ and
$\varepsilon
=(-1)^{n+1}$;
\end{itemize}

  \item[\rm(b)]
of the following
matrix pairs that are
given by
$0\ne\lambda\in\mathbb
F$ determined up to
replacement by
$\bar\lambda^{-1}$:
\begin{itemize}
  \item[\rm(i)]
$(J_n(\lambda)\oplus
J_n(\lambda)^{-*},
I_n\dia
 I_n)$ if
$|\lambda|\ne 1$,

  \item[\rm(ii)]
$(\lambda\Lambda_n,
\pm i^{n-1}F_n)$ if
$|\lambda|=1$;
\end{itemize}

\item[\rm(c)] of the
following matrix pairs
that are given by
$0\ne\lambda\in
\mathbb P+\mathbb Pi$
determined up to
replacement by
$\lambda^{-1}$ $($by
$\lambda^{-1}$,
$\bar\lambda,$ and
$\bar\lambda^{-1}$ in
{\rm(iii))}:
\begin{itemize}
  \item[\rm(i)]
$(J_n(\lambda)\oplus
J_n(\lambda)^{-T},
I_n\dia \varepsilon
I_n)$ if
$\lambda\in\mathbb P$,
except for
$\lambda=\pm 1$ and
$\varepsilon
=(-1)^{n+1}$,

  \item[\rm(ii)]
$(\lambda\Lambda_n,\pm
F_n)$ if $\lambda=\pm 1$%
\footnote{This gives
$4$ pairs:
$(\Lambda_n,F_n)$,
$(\Lambda_n,-F_n)$,
$(-\Lambda_n,F_n)$,
and
$(-\Lambda_n,-F_n)$.}
 and
$\varepsilon
=(-1)^{n+1}$,

\item[\rm(iii)]
$(J_n(\lambda)^{\mathbb
P}\oplus
(J_n(\lambda)^{\mathbb
P})^{-T}, I_{2n}\dia
\varepsilon I_{2n})$
if $\lambda
\notin\mathbb P$ and
$|\lambda|\ne 1$,

  \item[\rm(iv)]
$((\lambda\Lambda_n)
^{\mathbb P},\pm
(i^{n-(\varepsilon
+1)/2}F_n)^{\mathbb
P})$ if $\lambda
\notin\mathbb P$ and
$|\lambda|=1$;
\end{itemize}

 \item[\rm(d)]
of the following
matrix pairs that are
given by
$0\ne\lambda\in
\mathbb P+\mathbb Pi$
determined up to
replacement by
$\lambda^{-1}$,
$\bar\lambda,$ and
$\bar\lambda^{-1}$:
\begin{itemize}
  \item[\rm(i)]
$(J_n(\lambda)\oplus
J_n(\lambda)^{-*},
I_n\dia\varepsilon
I_n)$ if $|\lambda|\ne
1$,

  \item[\rm(ii)]
$(\lambda\Lambda_n,
\delta
i^{n-(\varepsilon
+1)/2}F_n)$ if
$|\lambda|=1$, where
\[
\delta:=
  \begin{cases}
     1,&
 \text{if $\lambda=
\pm 1$,
 the involution is
\eqref{ne},
$\varepsilon=(-1)^{n}$,
}\\
& \text{and if
$\lambda= \pm 1$,
 the involution is
\eqref{nen},
$\varepsilon=(-1)^{n+1}$};\\
    \pm 1,&
 \text{otherwise.}
  \end{cases}
\]
\end{itemize}
\end{itemize}
\end{theorem}

In this theorem
``determined up to
replacement by'' means
that a matrix pair
reduces by
transformations
\eqref{aaa} to the
matrix pair obtained
by making the
indicated replacement
(i.e., they give the
same $(A,B)$ but in
different bases).

\begin{remark}\label{remm}
The matrix
$i^{n-(\varepsilon
+1)/2}F_n$ in (c)(iv)
and (d)(ii) can be
replaced by $F_n$ if
$\varepsilon
=(-1)^{n+1}$ and by
$iF_n$ if $\varepsilon
=(-1)^{n}$. The pairs
\begin{equation*}\label{dpl}
(\lambda\Lambda_n, \pm
i^{n-1}F_n),\quad
((\lambda\Lambda_n)
^{\mathbb P},\pm
(i^{n-(\varepsilon
+1)/2}F_n)^{\mathbb
P}),\quad
(\lambda\Lambda_n,
\delta
i^{n-(\varepsilon
+1)/2}F_n)
\end{equation*}
in (b)(ii), (c)(iv),
and (d)(ii) can be
replaced by
\begin{equation*}\label{dpl2}
(\lambda\Omega_n, \pm
E_n),\quad
((\lambda\Omega_n)^{\mathbb
P},\pm
(\sqrt{\varepsilon}E_n)^{\mathbb
P}),\quad
(\lambda\Omega_n,
\delta\sqrt{\varepsilon}
E_n),
\end{equation*}
where $\sqrt{-1}=i$
and
\[
\Omega_n:=\begin{bmatrix}
1&2i&2i^2&\ddots&2i^{n-1}
\\&1&2i&\ddots&\ddots
\\&&1&\ddots&2i^2
\\&&&\ddots&2i
\\0&&&&1
\end{bmatrix},\quad
E_n:=\begin{bmatrix}
0&&1
\\
&\ddd&\\
1&&0
\end{bmatrix}\quad
\text{($n$-by-$n$).}
\]
This remark follows
from the proof of
Theorem \ref{theor}
and from the
equalities
\[
S_n^{-1}\Lambda_n S_n=
\Omega_n,\qquad
S_n^{*}i^{n-1}F_n S_n=
E_n,
\]
where
$S_n:=\diag(1,i,i^2,i^3,
\dots,i^{n-1})$ (i.e.,
$(\Lambda_n,
i^{n-1}F_n)$ and
$(\Omega_n,E_n)$ gave
the same $(A,B)$ but
in different bases).
\end{remark}

Due to the following
lemma, we have the
right to consider only
the involutions
\eqref{ne} and
\eqref{nen} on
$\mathbb H$.

\begin{lemma}\label{l1a}
Let\/ $\mathbb H$ be
the skew field of
quaternions over a
real closed field\/
$\mathbb P$. If any
involution on\/
$\mathbb H$ is not
quaternionic
conjugation
\eqref{ne}, then it
becomes quaternionic
semiconjugation
\eqref{nen} after a
suitable reselection
of the imaginary units
$i,j,k$.\footnote{Remark at proofreading: this statement was proved in [Randow, The involutory
antiautomorphisms of the
quaternion algebra,
\emph{Amer. Math. Monthly} 74
(1967) 699--700].}
\end{lemma}

\begin{proof}
The {\it absolute
value} of a quaternion
$h=a+bi+cj+dk$ is the
unique nonnegative
real root
\[
|h|:=\sqrt{a^2+b^2+c^2+d^2}=\sqrt{h\bar
h}\in\mathbb P,
\]
where $\bar
h:=a-bi-cj-dk$ is the
\emph{conjugate
quaternion}
($a^2+b^2+c^2+d^2$ is
a square by Lemma
\ref{l00}(a)). Then
$h^{-1}=|h|^{-2}\bar
h$ if $h\ne 0$.

The vector space of
purely imaginary
quaternions
\[
{\mathbb E}:=
\{bi+cj+dk\,|\,b,c,d\in\mathbb
P\}
\]
can be considered as
the Euclidean space
over $\mathbb P$ with
scalar product
\[
(bi+cj+dk,
b'i+c'j+d'k):=
bb'+cc'+dd'.
\]
Then $\{i,j,k\}$ is an
orthonormal basis,
$|h|$ is the length of
$h\in{\mathbb E}$, and
the multiplication of
two purely imaginary
quaternions can be
represented in the
form
\begin{equation}\label{es}
h_1h_2=[h_1,h_2]-(h_1,h_2),
\qquad
h_1,h_2\in{\mathbb E},
\end{equation}
where $[h_1,h_2]$ is
the vector product (if
$\mathbb P=\mathbb R$
then we may use its
geometrical
definition, otherwise
we use its definition
via determinants) and
$(h_1,h_2)$ is the
scalar product; in
particular, $[i,j]=k$
and $(i,j)=0$.

If $\{i',j'\}$ is a
pair of orthonormal
quaternions in
${\mathbb E}$ (i.e.,
$|i'|=|j'|=1$ and
$(i',j')=0$), then
$i'$, $j'$, $k':=i'j'$
can be taken as a new
set of imaginary
units.

Let $h\mapsto \widehat
h$ be an involution on
$\mathbb H$ that is
different from
quaternionic
conjugation
\eqref{ne}. Let us
prove that it acts
identically on
$\mathbb P$; that is,
$\widehat r=r$ for all
$r\in\mathbb P$. Each
$r\in\mathbb P$
commutes with all
$h\in\mathbb H$, hence
$\widehat r$ commutes
with all $\widehat
h\in\mathbb H$. Since
$\mathbb P$ is the
center of $\mathbb H$,
$\widehat r\in\mathbb
P$ and $r\mapsto
\widehat r$ is an
involution on $\mathbb
P$. If $\mathbb
P_{\circ}:=\{r\in\mathbb
P\,|\, \widehat r=r\}$
is its fixed field,
then the algebraically
closed field $\mathbb
P+\mathbb Pi$ has a
finite degree over
$\mathbb P_{\circ}$,
and so $\mathbb
P_{\circ}$ is a real
closed field. By Lemma
\ref{l00}(a), this
degree is 2, and so
$\mathbb
P_{\circ}=\mathbb P$.

Since $h\mapsto
\widehat h$ is not
quaternionic
conjugation, by
\cite[Chapter 8, \S
11, Proposition
2]{bour2} there exists
$h=a+bi+cj+dk\notin
\mathbb P$ such that
$\widehat h=h$. Put
$$e:=(b^2+c^2+d^2)^{-1/2}(bi+cj+dk),$$
then $e\in{\mathbb
E}$, $|e|=1$, and
$\widehat e=e$.

Choose any
$f\in{\mathbb E}$ of
length $1$ being
orthogonal to $e$.
Then by \eqref{es}
\begin{equation}\label{ew}
 e^2=f^2=-1,\quad ef=-fe.
\end{equation}
Write $\mathbb
K:=\mathbb P+\mathbb P
e$.  Since
$\{e,f,ef\}$ is a
basis of ${\mathbb
E}$, there are
$a,b,c,d\in \mathbb P$
such that
\begin{equation}\label{ey}
 \widehat f=a+be+cf+def=\varepsilon
 +\delta f,\quad
 \varepsilon:=a+be,\:
 \delta:=c+de\in\mathbb
 K.
\end{equation}
Then $$
\widehat{\widehat f\
}=\widehat
{\varepsilon}+
\widehat f\widehat
{\delta}={\varepsilon}+
\widehat
f{\delta}={\varepsilon}+
(\varepsilon
 +\delta f){\delta}=
{\varepsilon}+
\varepsilon {\delta}+
 \delta f{\delta}.
$$
By \eqref{ew},
$fe=-ef$, and so
$\widehat{\widehat f\
}=
{\varepsilon}(1+\delta)+
 \delta {\delta}' f$
with
${\delta}':=c-de$. But
$\widehat{\widehat f\
}=f$, hence
$f={\varepsilon}(1+{\delta})
+\delta{\delta}'\! f$.
Since
$\varepsilon(1+{\delta}),
\delta{\delta}'\in\mathbb
K$, and $\mathbb
H=\mathbb K+\mathbb
Kf$, we have
\begin{equation*}\label{smpf}
\text{${\varepsilon}=0$
\ \ or \
${\delta}=-1$,\qquad
and\qquad
$\delta{\delta}'=c^2+d^2=1$.}
\end{equation*}
\medskip

{\it Case 1:
${\varepsilon}=0$}.
Then $\widehat
f=\delta f$. Since
$\mathbb K$ is the
algebraic closure of
$\mathbb P$, there
exist $x,y\in \mathbb
P$ such that
$(x+ye)^2=c+de=\delta$.
In view of $e^2=-1$,
\[
(x^2+y^2)^2=((x+ye)(x-ye))^2
=(c+de)(c-de)=c^2+d^2=1.
\]
Thus $x^2+y^2=1$. Let
us write $k':=(x+ye)f$
and prove that the
quaternions $i':=ek',\
j':=e,\ k'$ form a
desired set of
imaginary units.

It suffices to check
that they are purely
imaginary quaternions
satisfying
\begin{equation}\label{lf5}
 |e|=|k'|=1,\qquad (e, k')=0,
\end{equation}
and that the
involution $h\mapsto
\widehat h$ has the
form \eqref{nen} with
respect to these
imaginary units; i.e.,
\begin{equation}\label{msst}
\widehat{ek'}=-ek',\qquad
\widehat e=e,\qquad
\widehat{k'}=k'.
\end{equation}

By \eqref{ew}
\begin{align*}
k^{\prime\,
2}&=(x+ye)f(x+ye)f=
(x+ye)(x-ye)f^2\\&
=(x^2+y^2)f^2=f^2 =-1.
\end{align*}
In view of \eqref{es},
$(k',k')=1$, and so
$|k'|=1$. The
inclusion
\[
ek'=e(x+ye)f=xef-yf\in{\mathbb
E}
\]
implies $(e,k')=0$.
This proves
\eqref{lf5}.

Furthermore,
\begin{align*}
\widehat{k'}&=\widehat
f(x+ye)=(x+ye)^2f(x+ye)=
(x+ye)^2(x-ye)f\\&
=(x+ye)(x^2+y^2)f
=(x+ye)f=k'
\end{align*}
and
$\widehat{ek'}=k'e=-ek'$.
This proves
\eqref{msst}.
\medskip

{\it Case 2:
$\delta=-1$}. Let us
prove that the
quaternions
\begin{equation*}\label{ktfi}
i':=f, \quad
j':=e,\quad k':=fe
\end{equation*}
form a desired set of
imaginary units. The
conditions
$|f|=|e|=1$, $(f,
e)=1$, and $\widehat
e=e$ hold.

By \eqref{ey},
$\widehat
f=\varepsilon-f=a+be-f$.
In view of \eqref{ew},
${fe}=-{ef}$,
$\widehat{fe}=-\widehat{ef}$,
$e\widehat f=
-\widehat fe$, and so
$
e(a+be-f)=-(a+be-f)e$.
Since $-ef=fe$, we
have $(a+be)e=0$,
hence $\widehat f=-f$.
Finally,
$\widehat{fe}=\widehat
e\widehat f=-ef=fe$.
\end{proof}

\subsection{Isometric operators over a field
of characteristic
different from 2}
\label{sub_02}

Canonical matrices of
pairs $(A,B)$
consisting of a
nondegenerate
Hermitian or
skew-Hermitian form
$B$ and an isometric
operator $A$ were
obtained in
{\cite[Theorem
5]{ser_izv}} over any
field $\mathbb F$ of
characteristic
different from 2 up to
classification of
Hermitian forms. They
were based on the
Frobenius canonical
matrices for
similarity. We
rephrase \cite[Theorem
5]{ser_izv} in Theorem
\ref{Theorem 5} from
this section in terms
of an \emph{arbitrary}
set of canonical
matrices for
similarity. This
flexibility will be
used in the proof of
Theorem \ref{theor}.
An analogous
flexibility was used
in \cite{hor-ser} to
simplify over $\mathbb
C$ the canonical
matrices for
congruence and
*congruence from
\cite[Theorems
3]{ser_izv} (a direct
proof that the
matrices from
\cite{hor-ser} are
canonical is given in
\cite{hor-ser_regul,
hor-ser_can}).

For each polynomial
\[
f(x)=a_0x^n+a_1x^{n-1}+\dots
+a_n\in \mathbb F[x],
\]
we define the
polynomials
\begin{align*}\label{mau}
\bar f(x)&:=\bar
a_0x^n+\bar
a_1x^{n-1}+\dots+\bar
a_n,\\
f^{\vee}(x)&:=\bar
a_n^{-1}(\bar
a_nx^n+\dots+\bar
a_1x+\bar
a_0)\quad\text{if }
a_n\ne 0.
\end{align*}

\begin{lemma}[{\cite[Lemma 6]{ser_izv}}]
\label{LEMMA 7} Let
$\mathbb F$ be a field
with involution
$a\mapsto \bar a$
$($possibly, the
identity$)$, let $p(x)
= p^{\vee}(x)$ be an
irreducible polynomial
over $\mathbb F$, and
let $r$ be the integer
part of $(\deg
p(x))/2$. Consider the
field
\begin{equation}\label{alft}
\mathbb F(\kappa) =
\mathbb
F[x]/p(x)\mathbb
F[x],\qquad \kappa:=
x+p(x)\mathbb F[x],
\end{equation}
with involution
\begin{equation}\label{alfta}
f(\kappa)^{\circ} :=
\bar f(\kappa^{-1}).
\end{equation}
Then each element of\/
$\mathbb F(\kappa)$ on
which the involution
acts identically is
uniquely representable
in the form
$q(\kappa)$, where
\begin{equation}\label{ser13}
q(x)=a_rx^r+\dots+
a_1x +a_0+\bar
a_1x^{-1}+\dots+\bar
a_rx^{-r},
    \quad a_0 = \bar a_0,
\end{equation}
$a_0,\dots
a_r\in\mathbb F,$ and
if $\deg p(x) = 2r$
then
\begin{equation}\label{uvp}
a_r=
  \begin{cases}
    0
&\text{if the
involution
on $\mathbb F$ is the
identity}, \\
    \bar a_r
&\text{if the
involution on $\mathbb
F$ is not the identity
and
$p(0)\ne 1$},\\
    -\bar a_r
&\text{if the
involution on $\mathbb
F$ is not the identity
and $p(0)=1$}.
  \end{cases}
\end{equation}
\end{lemma}

\begin{proof}
\emph{Case 1:} $\deg
p(x) = 2r +1$. The
elements
$\kappa^{r},\dots,1,\dots,
\kappa^{-r}$ ($\kappa$
is defined in
\eqref{alft}) form a
basis of $\mathbb
F(\kappa)$ over
$\mathbb F$.
Therefore, each
element of $\mathbb
F(\kappa)$ is uniquely
representable in the
form
\begin{equation}\label{df8i}
a_{r}\kappa^{r}+\dots
    +a_0+\dots+a_{-r}\kappa^{-r},
    \qquad a_r,\dots,
a_{-r}\in\mathbb F.
\end{equation}
The involution
\eqref{alfta} acts
identically on
\eqref{df8i} if and
only if $a_{i} = \bar
a_{-i}$ for all
$i=0,1,\dots,r$.
\medskip

\emph{Case 2: $\deg
p(x) = 2r$ and the
involution on $\mathbb
F$ is the identity.}
Then the involution
\eqref{alfta} acts
identically on the
elements
$$
a_{r-1}\kappa^{r-1}+\dots
+a_0+\dots
+a_{r-1}\kappa^{-r+1},
\qquad
a_0,\dots,a_{r-1}\in\mathbb
F;$$ they are distinct
and form over $\mathbb
F$ a subspace of
dimension $r$, which
is contained in the
fixed field
\begin{equation}\label{iiyd}
\mathbb
F(\kappa)_{\circ}:=\bigl\{f(\kappa)
\in{\mathbb
F(\kappa)}\,\bigr|\,
\overline{f(\kappa)}=f(\kappa)
\bigr\}
\end{equation}
of $\mathbb F(\kappa)$
with respect to the
involution
\eqref{alfta}.
$\mathbb
F(\kappa)_{\circ}$ has
the same dimension $r$
over $\mathbb F$
because $\dim_{\mathbb
F}\mathbb
F(\kappa)=2r$, and so
the subspace and the
fixed field coincide.

\emph{Case 3: $\deg
p(x) = 2r$ and the
involution on $\mathbb
F$ is not the
identity.} Let
\begin{equation*}\label{lucf}
p(x)=x^{2r}+p_1x^{2r-1}
+\dots+
p_{2r-1}x+p_{2r},
\end{equation*}
then
\begin{equation*}\label{lacf}
p^{\vee}(x)=\bar
p_{2r}^{-1}(\bar
p_{2r}x^{2r}+\bar
p_{2r-1}x^{2r-1}+\dots+
\bar p_{1}x+1).
\end{equation*}

The equality $p(0) =
p^{\vee}(0)$ implies
$p_{2r}=\bar
p_{2r}^{-1}$. Taking
any $b\in \mathbb F$
for which $\bar b\ne
b$ and putting
\[
\delta:=
  \begin{cases}
    1 +
\bar p_{2r} & \text{if
} p_{2r}\ne -1,
 \\
    b-\bar b& \text{if }
    p_{2r}=-1,
  \end{cases}
\]
we find that $\delta
p_{2r}=\bar\delta$.
Then $\delta \bar
p_{2r}^{-1}= \delta
p_{2r} =\bar\delta$,
and so
\[
\delta x^{-r}
p^{\vee}(x)=
\overline{\delta
p}_{2r}x^{r}
+\overline{\delta
p}_{2r-1}x^{r-1}
+\dots+
\overline{\delta
p}_{1}x^{1-r}
+\overline{\delta}x^{-r}.
\]
Since
\[
\delta x^{-r}
p^{\vee}(x)= \delta
x^{-r} p(x)=\delta
x^{r} +\delta
p_1x^{r-1}+\dots+
\delta
p_{2r-1}x^{1-r}+\delta
 p_{2r}x^{-r},
\]
the function $\pi(x):
= \delta x^{-r}p(x)$
has the form
\[
 \pi(x)=c_rx^r +\dots+c_1x
 +c_0+\bar c_{1}x^{-1}+\dots
 +\bar c_{r}x^{-r},
 \qquad   c_{0}=\bar
 c_0,\ c_r\ne 0.
\]
Using the equalities
$c_r=\delta$ and
$\delta
p_{2r}=\bar\delta$, we
find that $c_r
p_{2r}=\bar c_r$,
\begin{equation}\label{bvdsj}
c_r\ne
  \begin{cases}
    \bar c_r & \text{if
$p(0)=p_{2r}\ne 1$}, \\
    -\bar c_r &
    \text{if $p(0)=p_{2r}= 1$}.
  \end{cases}
\end{equation}

Let $q(x)$ be of the
form \eqref{ser13},
and let $q(\kappa) =
0$. Let us prove that
$q(x)=0$. We have
\[
\kappa^rq(\kappa) =
0,\quad x^rq(x) \equiv
0\bmod p(x),\quad
x^rq(x) = ap(x)
\]
for some $a\in\mathbb
F$. Thus,
\[
q(x) =
a\delta^{-1}\delta
x^{-r}p(x)=b
\pi(x),\qquad
b:=a\delta^{-1};
\]
equating the first
coefficients and
equating the last
coefficients, we
obtain $a_r= b c_r$
and $\bar a_r= b \bar
c_r$. So $b = \bar b$
and in view of
\eqref{uvp} and
\eqref{bvdsj} the
equality $q(x) =
b\pi(x)$ is possible
only if $q(x) = 0$.

Consequently, the
elements $q(\kappa)$
with $q(x)$ of the
form \eqref{ser13}
belong to
\eqref{iiyd}, they are
distinct and form a
vector space of
dimension $2r$ over
the fixed field
$\mathbb
F_{\circ}=\{a\in{\mathbb
F}\,|\, \bar{a}=a\}$
of $\mathbb F$. But
this is the dimension
over $\mathbb
F_{\circ}$ of the
whole fixed field
\eqref{iiyd}, so the
vector space coincides
with \eqref{iiyd}.
\end{proof}

Two $n\times n$
matrices $M$ and $N$
are said to be
\emph{similar} or
*\!\emph{congruent} if
$S^{-1}MS=N$ or
$S^{*}MS=N$,
respectively, for some
nonsingular $S$.

We say that a square
matrix is
\emph{indecomposable
for similarity} if it
is not similar to a
direct sum of square
matrices of smaller
sizes. Let ${\cal
O}_{\mathbb F}$ be any
maximal set of
nonsingular
indecomposable
canonical matrices for
similarity;\label{papage}
this means that each
nonsingular
indecomposable matrix
is similar to exactly
one matrix from ${\cal
O}_{\mathbb F}$.

For example, ${\cal
O}_{\mathbb F}$ may
consist of all
nonsingular {\it
Frobenius
blocks}---i.e., the
matrices
\begin{equation}\label{3.lfo}
\Phi=\begin{bmatrix}
0&&
0&-c_n\\1&\ddots&&\vdots
\\&\ddots&0&-c_2\\
0&&1& -c_1
\end{bmatrix}
\end{equation}
whose characteristic
polynomials
$\chi_{\Phi}(x)$ are
powers of irreducible
polynomials
$p_{\Phi}(x)\ne x$:
\begin{equation}\label{ser24}
\chi_{\Phi}(x)=p_{\Phi}(x)^s
=x^n+
c_1x^{n-1}+\dots+c_n.
\end{equation}
If $\mathbb F$ is an
algebraically closed
field, then ${\cal
O}_{\mathbb F}$ may
consist of all
nonsingular Jordan
blocks.

For $\varepsilon =\pm
1$ and each
nonsingular matrix
$\Phi$ that is
indecomposable for
similarity, if there
exists a nonsingular
$M$ satisfying $M
=\varepsilon
M^*=\Phi^*M\Phi$ then
we fix any and denote
it by
$\Phi_{(\varepsilon)}$
(we follow the
notation in
\cite{ser_izv}).

It suffices to
construct
$\Phi_{(\varepsilon)}$
only for the matrices
$\Phi\in{\cal
O}_{\mathbb F}$
because if
$\Phi_{(\varepsilon)}$
exists and $\Psi$ is
similar to $\Phi$ then
$\Psi_{(\varepsilon)}$
also exists: if
\begin{equation}\label{mshd}
\Phi_{(\varepsilon)}
=\varepsilon
\Phi_{(\varepsilon)}^*
=\Phi^*\Phi_{(\varepsilon)}\Phi,
\end{equation}
then we can take
\begin{equation}\label{ndw}
\Psi_{(\varepsilon)}
=S^*\Phi_{(\varepsilon)}S,\qquad
 \Psi =S^{-1}\Phi S
\end{equation}
and obtain
\begin{equation}\label{mskhd}
\Psi_{(\varepsilon)}
=\varepsilon
\Psi_{(\varepsilon)}^*
=\Psi^*\Psi_{(\varepsilon)}\Psi.
\end{equation}
Moreover, if
$\Psi_{(\varepsilon)}$
is any matrix that is
*congruent to
$\Phi_{(\varepsilon)}$,
then it satisfies
\eqref{mskhd} with
$\Psi$ defined by
\eqref{ndw}.

Existence conditions
and an explicit form
of
$\Phi_{(\varepsilon)}$
for Frobenius blocks
$\Phi$ over a field of
characteristic not $2$
were established in
Theorem 9 of
\cite{ser_izv}; this
result is represented
in Lemma \ref{lsdy1}
with a detailed proof.
Over algebraically or
real closed fields, we
construct in Lemma
\ref{lsdy} matrices
$\Psi_{(\varepsilon)}$
that are *congruent to
$\Phi_{(\varepsilon)}$
from Lemma \ref{lsdy1}
but are much simpler.

Theorem 5 of
{\cite{ser_izv}},
which was formulated
only for the set
${\cal O}_{\mathbb F}$
of all Frobenius
blocks, is extended to
any ${\cal O}_{\mathbb
F}$ in the following
theorem.

\begin{theorem}
\label{Theorem 5} Let
$A$ be an isometric
operator on a
finite-dimensional
vector space with
nondegenerate
$\varepsilon$-Hermitian
form $B$ over a field
$\mathbb F$ of
characteristic
different from $2$.
Let ${\cal O}_{\mathbb
F}$ be a maximal set
of nonsingular
indecomposable
canonical matrices for
similarity over
$\mathbb F$.
 Then the pair $(A, B)$
can be given in some
basis by a direct sum
of matrix pairs of the
following types:
\begin{itemize}
  \item[\rm(i)]
$(\Phi\oplus\Phi^{-*},
I\dia \varepsilon I)$,
where $\Phi\in{\cal
O}_{\mathbb F}$ is
such that
$\Phi_{(\varepsilon)}$
does not exist $($see
Lemma
{\rm\ref{lsdy1}(a))}.

  \item[\rm(ii)]
${\cal
A}^{q(x)}_{\Phi}:=(\Phi,
\Phi_{(\varepsilon)}q(\Phi))$,
where $\Phi\in{\cal
O}_{\mathbb F}$ is
such that
$\Phi_{(\varepsilon)}$
exists and $q(x)\ne 0$
is of the form
\eqref{ser13} in which
$r$ is the integer
part of $(\deg
p_{\Phi}(x))/2$. Here
$p_{\Phi}(x)$ is the
irreducible divisor of
the characteristic
polynomial of $\Phi$.
\end{itemize}
The summands are
determined to the
following extent:
\begin{description}
  \item [Type (i)]
up to replacement of
$\Phi$ by
$\Psi\in{\cal
O}_{\mathbb F}$ that
is similar to
$\Phi^{-*}$ $($i.e.,
whose characteristic
polynomial is
$\chi_{\Psi}(x)
=\chi_{\Phi}^{\vee}(x))$.

  \item [Type (ii)]
up to replacement of
the whole group of
summands
\[
{\cal
A}_{\Phi}^{q_1(x)}
\oplus\dots\oplus
{\cal
A}_{\Phi}^{q_s(x)}
\]
with the same $\Phi$
by
\[
{\cal
A}_{\Phi}^{q'_1(x)}
\oplus\dots\oplus
  {\cal A}_{\Phi}^{q'_s(x)}
\]
such that each
$q'_i(x)$ is a nonzero
function of the form
\eqref{ser13} and the
Hermitian forms
\begin{gather*}
q_1(\kappa)x_1^{\circ}x_1+\dots+
q_s(\kappa)x_s^{\circ}x_s,
\\
q'_1(\kappa)x_1^{\circ}x_1+\dots+
q'_s(\kappa)x_s^{\circ}x_s
\end{gather*}
are equivalent over
the field \eqref{alft}
with involution
\eqref{alfta}.
\end{description}
\end{theorem}

The proof of this
theorem given in
Section \ref{s_pro} is
a light modification
of the proof of
Theorem 5 in
{\cite{ser_izv}}.

Let $$f(x)=
\gamma_0x^m +
\gamma_1x^{m-1}+\dots+\gamma_m\in
\mathbb F[x],\qquad
\gamma_0\ne
0\ne\gamma_m.$$ A
vector $(a_1,
a_2,\dots, a_n)$ over
$\mathbb F$ is said to
be
\emph{$f$-recurrent}
if $n\le m$, or if
\[
\gamma_0 a_{l} +
\gamma_{1}a_{l+1}+\dots+
\gamma_ma_{l+m}=0,\qquad
l=1,2,\dots,n - m
\]
(by definition, it is
not $f$-recurrent if
$m=0$). Thus, this
vector is completely
determined by any
fragment of length
$m$.

The following lemma
was proved sketchily
in \cite[Theorem
9]{ser_izv}.

\begin{lemma}
\label{lsdy1}
 Let $\mathbb F$ be
a field of
characteristic
different from $2$
with involution
$($possibly, the
identity$)$. Let a
matrix $\Phi\in\mathbb
F^{n\times n}$ be
nonsingular and
indecomposable for
similarity; thus, its
characteristic
polynomial is a power
of some irreducible
polynomial
$p_{\Phi}(x)$.

{\rm(a)}
$\Phi_{(\varepsilon)}$
exists if and only if
\begin{equation}\label{lbdr}
 p_{\Phi}(x) =
p_{\Phi}^{\vee}(x),\
\text{and}
\end{equation}
\begin{equation}\label{4.adlw}
\parbox{20em}
{if the involution on
$\mathbb F$ is the
identity\\ and
$\varepsilon=(-1)^n$,
then $\deg p_{\Phi}(x)
> 1$.}
\end{equation}

{\rm(b)} If
\eqref{lbdr} and
\eqref{4.adlw} are
satisfied and if
$\Phi$ is a
nonsingular Frobenius
block \eqref{3.lfo}
with characteristic
polynomial
\begin{equation}\label{ser24lk}
\chi_{\Phi}(x)=p_{\Phi}(x)^s
=x^n+
c_1x^{n-1}+\dots+c_n,
\end{equation}
then for
$\Phi_{(\varepsilon)}$
one can take the
Toeplitz matrix
\begin{equation}\label{okjd}
\Phi_{(\varepsilon)}:=
[a_{i-j}]=
\begin{bmatrix}
a_0
&a_{-1}&\ddots&a_{1-n}
\\a_{1}&a_0
&\ddots&\ddots
\\\ddots&\ddots&
\ddots&a_{-1}
\\ a_{n-1}&\ddots
&a_{1}&a_0
\end{bmatrix}
\end{equation}
whose vector of
entries
$(a_{1-n},a_{2-n},\dots,a_{n-1})$
is the
$\chi_{\Phi}$-recurrent
extension of the
vector
$v=(a_{-m},\dots,a_{m})$
of length
\[
2m+1=
  \begin{cases}
    n & \text{if $n$ is odd}, \\
    n+1 & \text{if $n$ is
even,}
  \end{cases}
\]
defined as follows:

\begin{itemize}
  \item[\rm(i)]
$v:=(c_n-\varepsilon,0,
\dots,0,
\varepsilon\bar{c}_n-1
)$ if $n$ is even and
$c_n\ne\varepsilon$
$($see
\eqref{ser24lk}$)$;

  \item[\rm(ii)]
$v:=(c_1,-1,0,
\dots,0,-1,c_1)$
$(v:=(c_1,-2,c_1)$ for
$n=2)$ if $n$ is even,
$c_n=\varepsilon$, and
the involution on
$\mathbb F$ is the
identity;

  \item[\rm(iii)]
$v:=(a-\bar a, 0,\dots
, 0, \bar a-a)$
$($with any
$a\in\mathbb F$ such
that $\bar a\ne a)$ if
$n$ is even,
$c_n=\varepsilon$, the
involution is not the
identity, and also if
$n$ odd, $p_{\Phi}(x)
= x + c$, $c^{n-1}=-1$
$($then the involution
is not the
identity$)$.

  \item[\rm(iv)]
$v:=(1, 0,\dots , 0,
\varepsilon)$ if $n$
is odd and
$p_{\Phi}(x) \ne x +
c$, $c^{n-1}=-1$.
\end{itemize}
\end{lemma}

\begin{proof}
(a) Let
$\Phi\in\mathbb
F^{n\times n}$ be
nonsingular and
indecomposable for
similarity. Let us
prove that if
$\Phi_{(\varepsilon)}$
exists then the
conditions
\eqref{lbdr} and
\eqref{4.adlw} are
satisfied; we prove
the converse statement
in (b).

Let
$A:=\Phi_{(\varepsilon)}$
exist. By
\eqref{mshd}, $ A
=\varepsilon A^*=
\Phi^*A\Phi. $ Since
$A\Phi
A^{-1}=\Phi^{-*}$, we
have
\begin{align*}
\chi_{\Phi}(x) &=
\det(xI-\Phi^{-*})=
\det(xI-\bar\Phi^{-1})=
\det((-\bar\Phi^{-1})(I
-x\bar\Phi))=
\\&=\det(-\bar\Phi^{-1})\cdot
x^n\cdot \det(x^{-1}I
-\bar\Phi)=\chi_{\Phi}^{\vee}(x),
\end{align*}
where $n\times n$ is
the size of $\Phi$. In
the notation
\eqref{ser24},
$p_{\Phi}(x)^s =
p_{\Phi}^{\vee}(x)^s$,
which verify
\eqref{lbdr}.

To prove
\eqref{4.adlw},
suppose that the
involution on $\mathbb
F$ is the identity.

If $\varepsilon =-1$
then $A=-A^T$. Since
$A$ is skew-symmetric
and nonsingular, $n$
is even and so
$\varepsilon\ne(-1)^n$.

Let $\varepsilon =1$
and $\deg
p_{\Phi}(x)=1$. The
matrix $A$ is
symmetric and by
\eqref{lbdr}
$p_{\Phi}(x)=x\pm 1$.
Due to
\eqref{mshd}--\eqref{mskhd},
we may assume that
$\Phi=J_n(\pm 1)$.
Then $ A=J_n(\pm
1)^TAJ_n(\pm 1),$
  $J_n(\pm 1)^{-T}
A=AJ_n(\pm 1),$ and
\begin{equation}\label{smsi}
\begin{bmatrix}
0&&&0\\-1&0&&\\&\ddots&\ddots&
\\ *&&-1&0\end{bmatrix}A=
A\begin{bmatrix}
0&1&&0\\&0&\ddots&\\&&\ddots&1
\\ 0&&&0\end{bmatrix}.
\end{equation}
This implies that
\[
A=\begin{bmatrix} 0&&a_n\\
&\ddd&
\\
a_1&& *
\end{bmatrix}
\]
for some
$a_1,\dots,a_n$. Then
by \eqref{smsi}
\begin{equation*}\label{sme}
\begin{bmatrix} 0&&&0\\
&&\ddd&-a_n\\
&0& \ddd
\\
0&-a_2&&*
\end{bmatrix}=
\begin{bmatrix} 0&&&0\\
&&\ddd&a_{n-1}\\
&0& \ddd
\\
0&a_1&&*
\end{bmatrix}
\end{equation*}
and \[
(a_1,a_2,\dots,a_n)
=(a_1,-a_1,a_1,\dots,(-1)^{n-1}a_1).
\]
Since $A$ is
symmetric, $a_1=a_n$.
If $n$ is even, then
$a_1=a_n=-a_1$, and so
$a_1=0$, contrary to
the nonsingularity of
$A$. Hence, $n$ is odd
and
$\varepsilon\ne(-1)^n$.

(b) Let $\Phi$ be a
nonsingular Frobenius
block \eqref{3.lfo}
with characteristic
polynomial
\eqref{ser24lk}
satisfying
\eqref{lbdr} and
\eqref{4.adlw}. Write
\begin{equation}\label{ser25}
\mu_{\Phi}(x):=
p_{\Phi}(x)^{s-1}=x^t+
b_1x^{t-1}+\dots+b_t,\qquad
b_0:=1.
\end{equation}

Let
\begin{equation}\label{oyuyf}
(a_{1-n},\dots,a_{n-1})
\end{equation}
be any vector that is
$\chi_{\Phi}$-recurrent
but is not
$\mu_{\Phi}$-recurrent.
Consider the matrix
$A:=[a_{i-j}]$ of the
form \eqref{okjd}. By
\eqref{lbdr},
\begin{equation}\label{lyf}
\begin{aligned}
&\chi_{\Phi}(x) =x^n+
c_1x^{n-1}+\dots+c_{n-1}x+c_n
\\=&\chi_{\Phi}^{\vee}(x)=
\bar c_n^{-1}(\bar
c_nx^n+ \bar
c_{n-1}x^{n-1}
+\dots+\bar c_1x+1),
\end{aligned}
\end{equation}
and so the last row of
$\Phi^*$ is
\[
(-\bar c_n,\dots,
-\bar
c_1)=c_n^{-1}(-1,
-c_1,\dots,-c_{n-1}).
\]
Hence
\begin{equation}\label{ufc}
\Phi^*A\Phi=
\Phi^*[a_{i-j-1}]=
[a_{i-j}]=A
\end{equation}
($a_n$ is defined by
this equality).

Let us show that $A$
is nonsingular. If
$w:= (a_{n-1},\dots,
a_0)$ is the last row
of $A$, then
\begin{equation}\label{ldyf}
w\Phi^{n-1},\
w\Phi^{n-2},\dots, w
\end{equation}
are the rows of $A$.
Suppose, on the
contrary, that they
are linearly
dependent. Then
$wf(\Phi) = 0$ for
some nonzero
polynomial $f(x)$ of
degree less than $n$.
If $p_{\Phi}(x)^r$ is
the greatest common
divisor of $f(x)$ and
$\chi_{\Phi}(x)=p_{\Phi}(x)^s$,
then $r<s$ and
$$p_{\Phi}(x)^r=f(x)g(x)+
\chi_{\Phi}(x)h(x)\qquad
\text{for some
}g(x),h(x)\in\mathbb
F[x].$$ Since
$wf(\Phi) = 0$ and
$w\chi_{\Phi}(\Phi) =
0$, $wp_{\Phi}(\Phi)^r
= 0$. So
$w\mu_{\Phi}(\Phi) =
0$. Because
\eqref{ldyf} are the
rows of $A$, for each
$i=0,1,\dots, n-t-1$
we have
\begin{align*}
(0,\dots,0,&b_0,\dots,b_t,
\underbrace{0,\dots,0}_{i})A\\
&=b_0w\Phi^{i+t}+
b_1w\Phi^{i+t-1}+\dots
+b_tw\Phi^{i} =
w\mu_{\Phi}(\Phi)\Phi^i
=0\Phi^i=0.
\end{align*}
Hence,
$(a_{1-n},\dots,a_{n-1})$
is
$\mu_{\Phi}$-recurrent,
a contradiction.

What is left is to
show that the vector
$v=(a_{-m},\dots,a_{m})$
defined in (i)--(iv)
is
$\chi_{\Phi}$-recurrent
but is not
$\mu_{\Phi}$-recurrent
because this will
imply that its
$\chi_{\Phi}$-recurrent
extension
\eqref{oyuyf} defines
the nonsingular matrix
$A=[a_{i-j}]$
satisfying
\eqref{ufc}; since $v$
has the form
\[
(\varepsilon\bar
a_m,\dots,\varepsilon
\bar a_1,a_0,a_1\dots,
a_m),\qquad
\varepsilon \bar
a_0=a_0,
\]
we have that
$A=\varepsilon A^*$
and so $A$ can be
taken for
$\Phi_{(\varepsilon)}$.

(i$'$) The vector (i)
of length $n+1$ is not
$\mu_{\Phi}$-recurrent.
By \eqref{lyf},
$c_n=\bar c_n^{-1}$.
The vector (i) is
$\chi_{\Phi}$-recurrent
since
 $
c_n-\varepsilon+
c_n(\varepsilon \bar
c_n-1)=0$.

(ii$'$) Let $n$ be
even,
$c_n=\varepsilon$, and
let the involution on
$\mathbb F$ be the
identity. Then
\eqref{lyf} implies
$\chi_{\Phi}(1)
=c_n^{-1}\chi_{\Phi}(1)$.

If $\chi_{\Phi}(1)=0$
then
$p_{\Phi}(x)=x-1$.
Since $n$ is even,
\begin{equation}\label{axe}
\varepsilon
=c_n=1=(-1)^n,
\end{equation}
contrary to
$\eqref{4.adlw}$.

Hence
$\chi_{\Phi}(1)\ne 0$.
This gives $c_n=1$,
and so $c_1=c_{n-1}$
by \eqref{lyf}. The
vector (ii) is
$\chi_{\Phi}$-recurrent
because
$c_1-c_1-c_{n-1}+c_nc_1=0$.

In the same way,
$\mu_{\Phi}(x) =
\mu_{\Phi}^{\vee}(x)$
implies $\mu_{\Phi}(1)
=b_t^{-1}\mu_{\Phi}(1)$
and so $b_t=1$. In
view of \eqref{axe},
the condition
\eqref{4.adlw} ensures
$\deg p_{\Phi}(x) >
1$, thus $\deg
\mu_{\Phi}(x)=t\le n -
2$. The vector (ii) is
not
$\mu_{\Phi}$-recurrent
since if $n>2$ then
its fragment $(-1,0,
\dots,0,-1)$ of length
$n-1$ is not
$\mu_{\Phi}$-recurrent
and if $n=2$ then
$\mu_{\Phi}(x)$ is a
scalar.

(iii$'$) Let first $n$
be even,
$c_n=\varepsilon$, and
the involution be not
the identity. Then
$c_n=\varepsilon=1$,
so the vector (iii) of
length $n+1$ is
$\chi_{\Phi}$-recurrent
and  is not
$\mu_{\Phi}$-recurrent.

Let now $n$ be odd,
$p_{\Phi}(x) = x + c$,
and $c^{n-1}= -1$.
Then the involution is
not the identity:
otherwise $p_{\Phi}(x)
= p_{\Phi}^{\vee}(x) =
x \pm 1$ contradicts
$c^{n-1}= -1$. The
vector (iii) is
$\chi_{\Phi}$-recurrent
because of its length
$n<n+1$. It is not
$\mu_{\Phi}$-recurrent
since
$\mu_{\Phi}(x)=(x+c)^{n-1}$,
and so
$b_t=c^{n-1}=-1$ in
\eqref{ser25}.

(iv$'$) Let $n$ be
odd, and if
$p_{\Phi}(x) = x + c$
then $c^{n-1}\ne -1$.
The vector (iv) is
$\chi_{\Phi}$-recurrent
since its length
$n<n+1$.

If $\deg p_{\Phi}(x)
> 1$ then the length
of the vector (iv) is
greater than $\deg
\mu_{\Phi}(x)=t+1$,
thus (iv) is not
$\mu_{\Phi}$-recurrent.

If $p_{\Phi}(x)=x+c$
then $b_t=c^{n-1}\ne
-1$. By
\eqref{4.adlw},
$\varepsilon =1$,
hence (iv) is not
$\mu_{\Phi}$-recurrent.
\end{proof}

\section{Systems of
forms and linear
mappings}\label{s_pos}

In this section we
present in detail the
method of articles
\cite{roi,ser_first,ser_izv}
for reducing the
problem of classifying
systems of forms and
linear mappings to the
problem of classifying
systems of linear
mappings.

Let $V$ be a vector
space over $\mathbb
F$. A mapping
$\varphi\colon V\to
\mathbb F$ is called
\emph{semilinear} if
$$
\varphi(ua+vb)=\bar{a}\varphi
(u)+ \bar{b}\varphi
(v)\qquad\text{for all
}\ u,v\in V,\ \ a,b\in
\mathbb F.
$$
The set of all
semilinear mappings on
$V$ is a vector space,
we call it the {\it
*dual space} to $V$
and denote by $V^*$.

We identify $V$ with
$V^{**}$ by
identifying $v\in V$
with $\varphi\mapsto
\overline{\varphi v}$,
$\varphi\in V^*.$

For every linear
mapping $A:U\to V$, we
define the {\it
*adjoint mapping}
$A^{*}\colon V^{*}\to
U^{*}$, in which $
A^{*}\varphi:=\varphi
A$ for all $\varphi\in
V^*.$

\subsection{Representations
of
dographs}\label{sub_pos1}

Classification
problems for systems
of linear mappings can
be formulated in terms
of quivers and their
representations
introduced by Gabriel
\cite{gab}. A
\emph{quiver} is an
oriented graph. Its
\emph{representation}
is given by assigning
to every vertex a
vector space and to
every arrow a linear
mapping of the
corresponding vector
spaces. To include
into consideration
systems of forms and
linear mappings, I
extended in
\cite{ser_first} the
notion of quiver
representations as
follows. A
\emph{dograph} (a
doubly oriented graph,
or an \emph{oriented
schema} in terms of
\cite{ser_izv}) is, by
definition, a graph
with nonoriented,
oriented, and doubly
oriented edges; for
example,
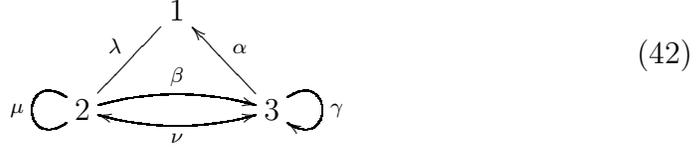
\begin{equation}\label{2.6}
\raisebox{20pt}{\xymatrix{
 &{1}&\\
 {2}\ar@(ul,dl)@{-}_{\mu}
 \ar@{-}[ur]^{\lambda}
  \ar@/^/@{->}[rr]^{\beta}
 \ar@/_/@{<->}[rr]_{\nu} &&{3}
 \ar[ul]_{\alpha}
 \ar@(ur,dr)^{\gamma}
 }}
\end{equation}

We suppose that the
vertices of each
dograph are
$1,2,\dots,n$, and
that there can be any
number of edges
between two vertices.

A {\it representation}
${\cal A}$ of a
dograph $D$ over
$\mathbb F$ is given
by assigning
\begin{itemize}
  \item
a vector space $V_i$
over $\mathbb F$ to
each vertex $i$,
  \item
a linear mapping
${A}_{\alpha}\colon
V_i\to V_j$ to each
arrow $\alpha\colon
i\to j$,

  \item
a sesquilinear form
${B}_{\beta}\colon
V_i\times V_j\to
{\mathbb F}$ to each
nonoriented edge
$\beta\colon i\lin\,
j\ (i\le
j)$\footnote{Thus,
$B_\beta$ is
semilinear on $V_i$
and linear on $V_j$ if
$i\le j$. This
condition is imposed
for definiteness and
it is unessential
because each
sesquilinear form
$B\colon U\times V\to
\mathbb F$ defines in
one-to-one manner the
sesquilinear form
$B^*\colon V\times
U\to \mathbb F$ as
follows: $B^*(v,u):=
\overline{B(u,v)}$.},
and

  \item
a sesquilinear form
${C}_{\gamma}\colon
V_i^*\times V_j^*\to
{\mathbb F}$ on the
*dual vector spaces
$V_i^*$ and $V_j^*$ to
each doubly oriented
edge $\gamma\colon
i\longleftrightarrow
j$ ($i\le j$).
\end{itemize}
Instead of $V_i$,
${A}_{\alpha}$,
${B}_{\beta}$,
$C_{\gamma}$ we
sometimes write ${\cal
A}_i$, ${\cal
A}_{\alpha}$, ${\cal
A}_{\beta}$, ${\cal
A}_{\gamma}$. The
\emph{dimension} of a
representation $\cal
A$ is the vector
\begin{equation}\label{bst}
\dim  {\cal A}:= (\dim
V_1,\dots,\dim V_n).
\end{equation}

For example, each
representation of the
dograph \eqref{2.6} is
a system
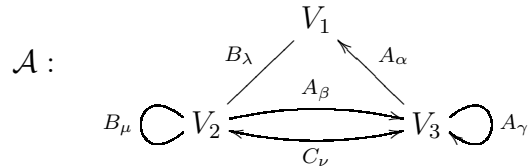
\begin{equation*}\label{2.6aa}
{\cal
A}:\quad\raisebox{20pt}{\xymatrix{
 &{V_1}&\\
 {V_2}
\save !<-2mm,0cm>
\ar@(ul,dl)@{-}_{B_{\mu}}
\restore
 \ar@{-}[ur]^{B_{\lambda}}
 \ar@/^/[rr]^{A_{\beta}}
 \ar@/_/@{<->}[rr]_{C_{\nu}} &&{V_3}
 \ar[ul]_{A_{\alpha}}
 \save
!<2mm,0cm>
\ar@(ur,dr)^{A_{\gamma}}
\restore }}
\end{equation*}
of vector spaces
$V_1,V_2,V_3$ over
$\mathbb F$, linear
mappings $
A_{\alpha}$,
$A_{\beta}$,
$A_{\gamma}$, and
forms
$$
B_{\lambda}\colon
V_1\times V_2\to
{\mathbb F},
 \quad
B_{\mu}\colon
V_2\times V_2\to
 {\mathbb F},\quad
 C_{\nu}\colon V^*_2\times
 V^*_3\to {\mathbb
 F}.
$$

A {\it morphism}
\begin{equation}\label{1.00}
f=(f_1,\dots,f_n):\:
{\cal A}\to{\cal A}'
\end{equation}
of representations
${\cal A}$ and ${\cal
A}'$ of $D$ is a set
of linear mappings
$f_i\colon V_i\to
V'_i$ that transform
$\cal A$ to ${\cal
A}'$; this means that
$$
f_j{A}_{\alpha}=
{A}'_{\alpha}{f}_i,\quad
 B_{\beta}(x,y)= B'_{\beta}
 (f_ix,f_jy),\quad
 C_{\gamma}(xf_i,yf_j)=
 C'_{\gamma}(x,y)
$$
for all oriented edges
$ \alpha\colon
i\longrightarrow j$,
nonoriented edges
$\beta\colon i\lin j\
(i\le j)$, and doubly
oriented edges
$\gamma\colon
i\longleftrightarrow
j\ (i\le j).$ The
composition of two
morphisms is a
morphism. A morphism
$f\colon{\cal
A}\to{\cal A}'$ is
called an {\it
isomorphism} and is
denoted by $f\colon
{\cal A}\is{\cal A}'$
if all $f_i\colon
V_i\to V'_i$  are
bijections. In this
case we say that
${\cal A}$ is
\emph{isomorphic} to
${\cal A}'$ and write
${\cal A}\simeq{\cal
A}'$. If ${\cal
A}={\cal A}'$, then
morphisms are called
\emph{endomorphisms}
and isomorphisms are
called
\emph{automorphisms}.

The {\it direct sum}\/
${\cal A}\oplus{\cal
A}'$ of
representations ${\cal
A}$ and ${\cal A}'$ of
$D$ is the
representation
consisting of the
vector spaces
$V_i\oplus V'_i$
(i=1,\,\dots,\,n), the
linear mappings
$$
{A}_{\alpha}\oplus
{A}'_{\alpha} :\
V_i\oplus V'_i\to
V_j\oplus V'_j,\qquad
\alpha:\,i\longrightarrow
j,
$$
and the forms
$$
B_{\beta}\oplus
B'_{\beta} :\
(V_i\oplus V'_i)\times
(V_j\oplus V'_j)\to
{\mathbb F},\qquad
\beta:\,i\lin\, j \
(i\le j),
$$
$$
{C}_{\gamma}\oplus
{C}'_{\gamma} :\
(V_i\oplus
V_i')^*\times
(V_j\oplus V_j')^*\to
{\mathbb F},\qquad
\gamma:\,i\longleftrightarrow\,
j \ (i\le j).
$$
A representation $\cal
A$ is {\it
indecomposable} if
$$
{\cal A}\simeq {\cal
B}\oplus{\cal
C}\quad\Longrightarrow\quad
{\cal B}=0 \ \text{ or
}\ {\cal C}=0,
$$
where $0$ is the
representation in
which all vector
spaces are $0$.

The set
$\rep(D,\mathbb F)$ of
representations of a
dograph $D$ over
$\mathbb F$ is a
category with
morphisms
\eqref{1.00}. But this
category is not
additive since the sum
of two morphisms
usually is not a
morphism. So the
properties of dograph
representations are
more complicated than
the properties of
quiver
representations, whose
morphisms form vector
spaces.

Let us denote by
$\Isom(D,\mathbb F)$
\label{jjjja} the
subcategory of
$\rep(D,\mathbb F)$
consisting of the same
objects and whose
morphisms are the
isomorphisms of
$\rep(D,\mathbb F)$.
Roiter \cite{roi}
proposed to study
representations of a
dograph $D$ embedding
$\Isom(D,\mathbb F)$
into the additive
category
$\rep(\underline{D},\mathbb
F)$ of representations
of some quiver
$\underline{D}$ with
involution. In Section
\ref{sub_pos2} we
introduce the notion
of a quiver with
involution and define
an involution on the
category of its
representations. In
Section \ref{sub_imb}
we construct the
embedding of
$\Isom(D,\mathbb F)$
to the category
$\rep(\underline{D},\mathbb
F)$. In Section
\ref{sub_relat} we
deal with dographs
with relations, they
admit to consider
systems of forms and
linear mappings
satisfying relations.
In Section
\ref{sub_pos3} we
reduce the problem of
classifying
representations of a
dograph $D$ with
relations to the
problems of
classifying
representations of the
quiver $\underline{D}$
with relations and
Hermitian forms over
finite extensions of
the center of $\mathbb
F$.

\subsection{Representations
of quivers with
involution}\label{sub_pos2}

By a \emph{quiver with
involution}, we mean a
quiver $Q$, in which
to every vertex $i$ we
associate some vertex
$i^*$ and to each
arrow $\alpha \colon
i\to j$ some arrow
$\alpha^* \colon
j^*\to i^*$ such that
$ i^*\ne i=i^{**}$ and
$\alpha^*\ne
\alpha=\alpha^{**}.$

The involution on $Q$
induces the following
involution on the
category of its
representations
$\rep(Q,\mathbb F)$:
\begin{itemize}
  \item
\emph{Involution on
representations.} To
each representation
$\cal M$ of $Q$ we
associate the
\emph{adjoint
representation} ${\cal
M}^{\circ}$ of $Q$
that assigns the
vector spaces ${\cal
M}^{\circ}_i:={\cal
M}^*_{i^*}$ and the
linear mappings ${\cal
M}^{\circ}_{\alpha}:={\cal
M}^*_{\alpha^*}$ to
all vertices $i$ and
arrows $\alpha$ of
$Q$.
  \item
\emph{Involution on
morphisms.} To each
morphism $f\colon
{\cal M}\to{\cal N}$
of representations of
$Q$ we associate the
{\it adjoint morphism}
\begin{equation}\label{kdtc}
f^{\circ}\colon {\cal
N}^{\circ}\to{\cal
M}^{\circ},\qquad
\text{in which }\
f^{\circ}_i:=f^*_{i^*}
\end{equation}
for all vertices $i$
of $Q$.
\end{itemize}
For example, consider
the quiver with
involution
\begin{equation*}\label{4.1jg}
\raisebox{23pt}{\xymatrix@R=4pt{
 &{2}\ar[dd]_{{\alpha}}
 \ar[ddrr]^(.25){{\beta}}&
 &{2^*} \\
 {Q:}&&\\
 &{1}\ar[uurr]^(.75){{\beta}^*}
 &&{1^*}\ar[uu]_{{\alpha}^*}
\ar@<-0.4ex>[ll]_{{\gamma}}
 \ar@<0.4ex>[ll]^{{\gamma}^*}
 }}
\end{equation*}
\begin{itemize}
  \item
For its representation
\begin{equation}\label{4.5}
\raisebox{23pt}{\xymatrix@R=4pt{
 &{U_1}\ar[dd]_{A_1}
 \ar[ddrr]^(.25){B_1}&
 &{U_2} \\
 {{\cal M}:}&&\\
 &{V_1}\ar[uurr]^(.75){B_2}
 &&{V_2}\ar[uu]_{A_2}
 \ar@<0.4ex>[ll]^{C_2}
 \ar@<-0.4ex>[ll]_{C_1}
 }}
\end{equation}
the adjoint
representation
\begin{equation*}\label{4.7}
\raisebox{23pt}{\xymatrix@R=4pt{
 &{U_2^*}\ar[dd]_{A_2^*}
 \ar[ddrr]^(.25){B_2^*}&
 &{U_1^*} \\
 {{\cal M}^{\circ}:}&&\\
 &{V_2^*}
 \ar[uurr]^(.75){B_1^*}
 &&{V_1^*}\ar[uu]_{A_1^*}
 \ar@<0.4ex>[ll]^{C_1^*}
 \ar@<-0.4ex>[ll]_{C_2^*}
 }}
\end{equation*}
is constructed as
follows: we replace
all vector spaces of
${\cal M}$ by the
*dual spaces, all
linear mappings by the
*adjoint mappings,
which reverses the
direction of each
arrow:
\begin{equation*}\label{4.7z}
\raisebox{23pt}{\xymatrix@R=4pt{
 &{U_1^*}\ar[dd];[]^{A_1^*}
 \ar[ddrr];[]_(.75){B_1^*}&
 &{U_2^*} \\
 {{\cal M}^{*}:}&&\\
 &{V_1^*}
 \ar[uurr];[]_(.25){B_2^*}
 &&{V_2^*}\ar[uu];[]^{A_2^*}
 \ar@<0.4ex>[ll];[]^{C_1^*}
 \ar@<-0.4ex>[ll];[]_{C_2^*}
 }}
\end{equation*}
rotate the obtained
representation around
the vertical axis, and
interchange $C_1^*$
and $C_2^*$.

  \item
For a morphism
\begin{equation}\label{4.12}
\raisebox{57pt}{\xymatrix@R=4pt{
 &&{U_1}\ar[dddd]_(.75){f_{2}}
 \ar[ddl]_{A_1}
 \ar[ddrrr]^(.25){B_1}&&&
 &{U_2}\ar[dddd]^{f_{2^*}}
 \\
 {\ {\cal M}:}
 \ar[dddd]_{f}
 \\
 &{V_1}\ar[dddd]_{f_{1}}
 \ar@<0.4ex>[uurrrrr]^(.75){B_2}
 \ar@<0.4ex>[rrrr];[]^{C_2}
 \ar@<-0.4ex>[rrrr];[]_{C_1}
 &&&&{V_2}\ar[dddd]^(.25){f_{1^*}}
 \ar[uur]_{A_2}&
 \\
 \\
 &&{\hat U_1}\ar[ddl]_{\hat A_1}
 \ar[ddrrr]^(.25){\hat B_1}&&&
 &{\hat U_2}
 \\
 {\ {\cal N}:}
 \\
 &{\hat V_1}\ar@<0.4ex>[uurrrrr]^(.63)
 {\hat B_2}
 \ar@<0.4ex>[rrrr];[]^{\hat C_2}
 \ar@<-0.4ex>[rrrr];[]_{\hat C_1}
 &&&&{\hat V_2}\ar[uur]_{\hat A_2}
 &
 }}\end{equation}
of its representations
${\cal M}$ and ${\cal
N}$, the adjoint
morphism
\begin{equation*}\label{4.13}
\raisebox{57pt}{\xymatrix@R=4pt
{
 &&{U_2^{*}}
 \ar[ddl]_{A_2^{*}}
 \ar[ddrrr]^(.25){B_2^{*}}&&&
 &{U_1^{*}}
 \\
 {\ \ \ {\cal M}^{\circ}:}
 \\
 &{V_2^{*}}
 \ar@<0.4ex>[uurrrrr]^(.75){B_1^{*}}
 \ar@<0.4ex>[rrrr];[]^{C_1^{*}}
 \ar@<-0.4ex>[rrrr];[]_{C_2^{*}}
 &&&&{V_1^{*}}
 \ar[uur]_{A_1^{*}}&
 \\
 \\
 &&{\hat U_2}^*\ar[uuuu]^(.25)
 {f_{2^*}^{*}}
 \ar[ddl]_{\hat A_2^{*}}
 \ar[ddrrr]^(.25){\hat B_2^{*}}&&&
 &{\hat U_1^{*}}\ar[uuuu]_{f_{2}^{*}}
 \\
 {\ \ \ {\cal N}^{\circ}:}
 \ar[uuuu]^{f^{\circ}}
 \\
 &{\hat V_2^{*}}\ar[uuuu]^{f_{1^*}^{*}}
 \ar@<0.4ex>[uurrrrr]^(.63){\hat B_1^{*}}
 \ar@<0.4ex>[rrrr];[]^{\hat C_1^{*}}
 \ar@<-0.4ex>[rrrr];[]_{\hat C_2^{*}}
 &&&&{\hat V_1^{*}}\ar[uuuu]_(.75){f_{1}^{*}}
 \ar[uur]_{\hat A_1^{*}}&
 }}\end{equation*}
is obtained as
follows: we replace
all vector spaces in
\eqref{4.12} by the
*dual spaces, all
linear mappings by the
*adjoint mappings,
rotate around the
vertical axis, and
interchange $C_1^*$
with $C_2^*$ and $\hat
C_1^*$ with $\hat
C_2^*$.
\end{itemize}

An isomorphism
$f\colon {\cal
M}\is{\cal N}$ of
selfadjoint
representations ${\cal
M}={\cal M}^{\circ}$
and ${\cal N}={\cal
N}^{\circ}$ is called
a {\it congruence} if
$f^{\circ}=f^{-1}$.

\subsection{Representations of
dographs as
selfadjoint
representations of
quivers with
involution}
\label{sub_imb}

For every dograph
${D}$, we denote by
$\underline{D}$ the
quiver with involution
obtained from $D$ by
replacing
\begin{itemize}
  \item
each vertex $i$ of
${D}$ by the vertices
$i$ and $i^*$,
  \item
each arrow
$\alpha\colon  i\to j$
by the arrows
$\alpha\colon i\to j$
and $\alpha^*\colon
j^*\to i^*$,
  \item
each nonoriented edge
$\beta\colon i\lin\,
j\ (i\le j)$ by the
arrows $\beta\colon
j\to i^*$ and
$\beta^*\colon i\to
j^*$,
  \item
each doubly oriented
edge $\gamma\colon
i\longleftrightarrow
j$ ($i\le j$)  by the
arrows $\gamma\colon
j^*\to i$ and
$\gamma^*\colon i^*\to
j$.
\end{itemize}
We define $i^{**}:=i$
and
$\alpha^{**}:=\alpha$
for all vertices $i$
and arrows $\alpha $
of the quiver
$\underline{D}$. For
example,
\begin{equation}\label{4.1}
\raisebox{20pt}{\xymatrix@R=4pt{
 &{2}\ar[dd]_{\alpha}
 \ar@{-}@/^/[dd]^{\beta}\\
{D}\colon &
  \\
 &{1}
\save !<1.5mm,0cm>
\ar@(ur,dr)@{<->}^{\gamma}
\restore}}
\qquad\qquad\qquad
\raisebox{23pt}{\xymatrix@R=4pt{
 &{2}\ar[dd]_{{\alpha}}
 \ar[ddrr]^(.25){{\beta}}&
 &{2^*} \\
 {\underline{D}\colon }&&\\
 &{1}\ar[uurr]^(.75){{\beta}^*}
 &&{1^*}\ar[uu]_{{\alpha}^*}
\ar@<-0.4ex>[ll]_{{\gamma}}
 \ar@<0.4ex>[ll]^{{\gamma}^*}
 }}
\end{equation}

The embedding of
$\Isom(D,\mathbb F)$
into $\rep(\underline
D,\mathbb F)$ (see
page \pageref{jjjja})
is constructed as
follows:

\begin{itemize}
  \item
\emph{Embedding of
representations.} To
each representation
$\cal A$ of $D$ over
$\mathbb F$, we
associate the
\emph{selfadjoint
representation}
$\underline {\cal A}$
of $\underline {D}$
obtained from $\cal A$
by replacing
\begin{itemize}
  \item
each vector space $V$
of ${\cal A}$ by the
spaces $V$ and $V^*$
(=\,the *dual space of
all semilinear forms
$V\to \mathbb F$),
  \item
each linear mapping
$A\colon U\to V$ by
the mutually *adjoint
mappings $A\colon U\to
V$ and $A^*\colon
V^*\to U^*$,

  \item
each sesquilinear form
$B\colon U\times
V\to\mathbb F$ by the
mutually *adjoint
mappings
\begin{equation*}\label{mse}
B\colon v\in V\mapsto
B(?,v)\in U^*,\qquad
B^*\colon u\in
U\mapsto
\overline{B(u,?)}\in
V^*,
\end{equation*}

  \item
each sesquilinear form
$C\colon U^*\times
V^*\to\mathbb F$ by
the mutually *adjoint
mappings
\begin{equation*}\label{zmse}
C\colon v^*\in
V^*\mapsto C(?,v^*)\in
U^{**}=U,\quad
C^*\colon u^*\in
U^*\mapsto
\overline{C(u^*,?)}\in
V.
\end{equation*}
\end{itemize}
(We use the same
letter for a
sesquilinear form
$B\colon U\times
V\to\mathbb F$ and for
the corresponding
mapping $B\colon V\to
U^*$. They have the
same matrices in any
bases $\{u_i\}$ of
$U$, $\{v_i\}$ of $V$,
and in the *dual basis
$\{u_i^*\}$ of $U^*$
defined by
$u^{*}_i(u_j)=0$ if
$i\ne j$ and
$u^{*}_i(u_i)=1$.)

For example, for the
dograph and the quiver
\eqref{4.1}:
\begin{equation}\label{4.2}
\raisebox{23pt}{\xymatrix@R=4pt{
 &{U}\ar[dd]_{A}
 \ar@{-}@/^/[dd]^{B}\\
{\cal A}\colon &
  \\
 &{V}
\save !<2mm,0cm>
\ar@(ur,dr)@{<->}^{C}
\restore}}
   \qquad   \qquad   \qquad
\raisebox{23pt}{\xymatrix@R=4pt{
 &{U}\ar[dd]_{A}
 \ar[ddrr]^(.25){B}&
 &{U^*} \\
 {\underline{\cal A}\colon }&&\\
 &{V}\ar[uurr]^(.75){B^*}
 &&{V^*}\ar[uu]_{A^*}
 \ar@<0.4ex>[ll]^{C^*}
 \ar@<-0.4ex>[ll]_{C}
 }}
\end{equation}

  \item
\emph{Embedding of
isomorphisms.} To each
isomorphism $f\colon
{\cal A}\is{\cal B}$
of representations of
a dograph ${D}$, we
associate the
\emph{congruence}
$\underline{f}\colon
\underline{\cal A}\is
\underline{\cal B}$ of
the corresponding
selfadjoint
representations of
$\underline{D}$ by
defining
$\underline{f}_{\,i}:=f_i$
and
$\underline{f}_{\,i^*}:=f_i^{-*}$
for each vertex $i$ of
$D$. For example, an
isomorphism
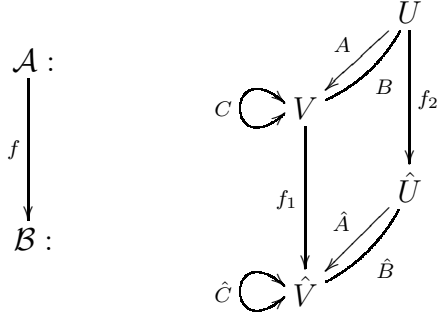
\begin{equation*}\label{4.11a}
\raisebox{57pt}{\xymatrix@R=4pt{
 &&&&{U}\ar[dddd]^{f_2}
 \ar[ddl]_{A}
 \ar@{-}@/^/[ddl]^{B}
  \\
{\ \cal
A}\ar[dddd]_{f}:&&&&
  \\
 &&&{V}\ar[dddd]_{f_1}
\save !<-2mm,0cm>
\ar@(ul,dl)@{<->}_{C}
\restore&
 \\
 \\
 &&&&{\hat U}\ar[ddl]_{\hat A}
 \ar@{-}@/^/[ddl]^{\hat B}
  \\
{\ \cal B}:&&&&
  \\
 &&&{\hat V}
\save !<-2mm,0cm>
\ar@(ul,dl)@{<->}_{\hat
C}& \restore
 }}
\end{equation*}
defines the congruence
\begin{equation*}\label{4.11b}
\raisebox{57pt}
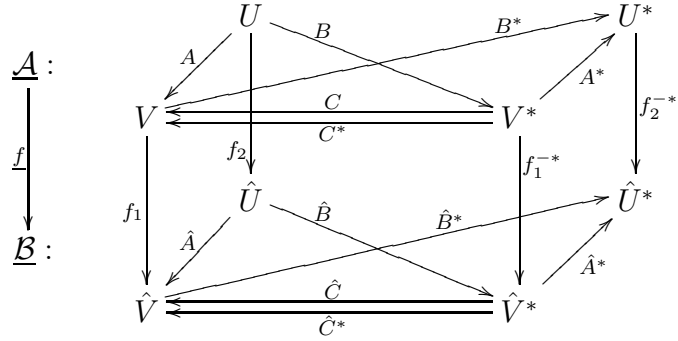
{\xymatrix@L=1pt@R=4pt{
 &&{U}\ar[dddd]_(.75){f_2}
 \ar[ddl]_{A}
 \ar[ddrrr]^(.25){B}&&&
 &{U^{*}}
 \ar[dddd]^{f_2^{-*}}
 \\
 {\ \underline{\cal A}:}
 \ar[dddd]_{\underline{f}}
 \\
 &{V}\ar[dddd]_{f_1}
 \ar@<0.4ex>[uurrrrr]^(.75){B^{*}}
 \ar@<0.4ex>[rrrr];[]^{C^*}
 \ar@<-0.4ex>[rrrr];[]_{C}
 &&&&{V^{*}}
 \ar[dddd]^(.25){f_1^{-*}}
 \ar[uur]_{A^{*}}&
 \\
 \\
 &&{\hat U}\ar[ddl]_{\hat A}
 \ar[ddrrr]^(.25){\hat B}&&&
 &{\hat U^{*}}
 \\
 {\ \underline{\cal B}:}
 \\
 &{\hat V}
 \ar@<0.4ex>[uurrrrr]^(.63){\hat B^{*}}
 \ar@<0.4ex>[rrrr];[]^{\hat C^*}
 \ar@<-0.4ex>[rrrr];[]_{\hat C}
 &&&&{\hat V^{*}}
 \ar[uur]_{\hat A^{*}}&
 }}
 \end{equation*}
\end{itemize}

Clearly, each
selfadjoint
representation of
$\underline{D}$ has
the form
$\underline{\cal A}$
and each congruence of
selfadjoint
representations has
the form
$\underline{f}\colon
\underline{\cal A}\to
\underline{\cal B}$.
Two representations
$\cal A$ and $\cal B$
of $D$ are isomorphic
if and only if the
corresponding
selfadjoint
representations
$\underline{\cal A}$
and $\underline{\cal
B}$ of $\underline{D}$
are congruent.
Therefore, \emph{the
problem of classifying
representations of a
dograph $D$ up to
isomorphism reduces to
the problem of
classifying
selfadjoint
representations of the
quiver $\underline{D}$
up to congruence.}

\subsection{Dographs with
relations}
\label{sub_relat}

A \emph{relation} on a
quiver $Q$ over a
field or skew field
$\mathbb F$ is a
formal expression of
the form
\begin{equation}\label{5.a1a}
\sum_{i=1}^m
c_i\alpha_{ip_i}\cdots
\alpha_{i2}\alpha_{i1}
=0,
\end{equation}
in which all $c_i$ are
nonzero elements of
the center of $\mathbb
F$ and
\begin{equation*}\label{alse}
 u\xrightarrow[\phantom{\qquad}]
 {\alpha_{i1}}
  u_{i2}
  \xrightarrow[\phantom{\qquad}]
  {\alpha_{i2}}
         \cdots
  \xrightarrow[\phantom{\qquad}]
  {\alpha_{i,p_i-1}}
  u_{ip_i} \xrightarrow[\phantom{\qquad}]
  {\alpha_{ip_i}}
  v
\end{equation*}
are oriented paths on
$Q$ with the same
initial vertex $u$ and
the same final vertex
$v$ ($u_{ij}$ and
$\alpha_{ij}$ are
vertices and arrows).
A path may have length
0 if $u=v$. This
``lazy'' path (without
arrows) is replaced by
$1$ in \eqref{5.a1a}
and gives a summand of
the form $c_i1$.
Therefore, if $u=v$
then \eqref{5.a1a} may
have `1' instead of
`0' in its right-hand
side.

A representation
${\cal A}$ of $Q$
\emph{satisfies the
relation
\eqref{5.a1a}} if
\begin{equation*}\label{5.a1au}
\sum_{i=1}^m c_i{\cal
A}_{\alpha_{ip_i}}\cdots
{\cal A}_{\alpha_{i2}}
{\cal A}_{\alpha_{i1}}
=0.
\end{equation*}
For example, the
problem of classifying
representations of the
quiver with relations
\[
\xymatrix{
 {1} \ar@(ul,dl)@{->}_{\alpha}
 \ar@(ur,dr)@{->}^{\beta}}\qquad\qquad
 \alpha \beta =\beta \alpha =0
\]
is the problem of
classifying pairs of
mutually annihilating
linear operators,
which was solved over
a field in
\cite{naz_bon}. The
notion of a quiver
with relations arose
in the theory of
representations of
finite dimensional
algebras over a field:
every algebra can be
given by a quiver with
relations and there is
a natural one-to-one
correspondence between
representations of the
algebra and
representations of the
quiver with relations.

By a {\it dograph with
relations}, we mean a
dograph $D$ with a
finite set of
relations on its
quiver with involution
$\underline D$, and
consider only those
representations $\cal
A$ of $D$, for which
the corresponding
selfadjoint
representations
$\underline{\cal A}$
of $\underline D$
satisfy these
relations. Clearly, if
$\underline{\cal A}$
satisfies the relation
\eqref{5.a1a}, then it
satisfies also the
\emph{adjoint
relation}
\begin{equation}\label{5.a1a3}
\sum_{i=1}^m
\bar{c}_i\alpha_{i1}^*
\alpha_{i2}^*\cdots
\alpha_{ip_i}^* =0.
\end{equation}

For example, the
problems of
classifying
representations of the
dographs
\begin{align}\nonumber
\xymatrix{
 {1}
 \ar@(ur,dr)@{-}^{\alpha}
 }&
          \\[2mm] \nonumber
\xymatrix{
 {1} \ar@(ul,dl)@{-}_{\alpha}
 \ar@(ur,dr)@{-}^{\beta}}&\qquad
 {\alpha=\varepsilon \alpha^*,
 \quad \beta=\delta \beta^*,}
           \\\label{2.8akda}
\xymatrix{
 {1\,} \ar@(ul,dl)@{->}_{\alpha}
 \ar@<-0.4ex>@(u,r)@{-}^{\beta}
\ar@<-0.4ex>@(r,d)@{<->}^{\gamma}
}\ \: &\qquad
  {\begin{matrix}
 \beta
 =\varepsilon \beta^*
 =\alpha^*\beta\alpha,\\
 \gamma\beta=1,\quad
 \beta\gamma=1,
 \end{matrix}}
               \\ \nonumber
 \xymatrix{
 {1\,} \ar@(ul,dl)@{->}_{\alpha}
 \ar@<-0.4ex>@(u,r)@{-}^{\beta}
\ar@<-0.4ex>@(r,d)@{<->}^{\gamma}
}\ \: &\qquad
  {\begin{matrix}
 \beta=\varepsilon \beta^*,\quad
  \beta\alpha=\alpha^*\beta,\\
 \gamma\beta=1,\quad
 \beta\gamma=1,
 \end{matrix}}
\end{align}
in which
$\varepsilon,\delta\in\{-1,1\}$
(due to the edges
$\gamma$ and the
relations
$\gamma\beta=1,$
$\beta\gamma=1$, the
form assigned to
$\beta$ in each
representation is
nondegenerate) are the
problems of
classifying,
respectively:
\begin{itemize}
  \item
sesquilinear forms,

  \item
pairs of forms, in
which the first is
$\varepsilon$-Hermitian
and the second is
$\delta$-Hermitian,
  \item
isometric operators on
a space with
nondegenerate
$\varepsilon$-Hermitian
form, and
  \item
selfadjoint operators
on a space with
nondegenerate
$\varepsilon$-Hermitian
form (an operator $A$
is \emph{selfadjoint}
with respect to $B$ if
$B(Au,v)=B(u,Av)$ for
all $u$ and $v$).
\end{itemize}

These problems were
solved in
\cite{ser_prep} and in
\cite[Theorems
3--6]{ser_izv} over
any field of
characteristic
different from 2 up to
classification of
Hermitian forms over
its finite extensions.
An analogous
description of pairs
of subspaces in a
space with an
indefinite scalar
product was given in
\cite{ser_sub} by
reducing it to the
problem of classifying
representations of the
dograph
\[
{\xymatrix@R=0,5pt{
&*{\ci}\ar@{->}[dl]
\\*{\ci}
\ar@(ul,dl)@{-}_{\alpha
}&&&\alpha^*=\varepsilon
\alpha.
\\
&*{\ci}\ar@{->}[ul]} }
\]

\subsection{Reduction
theorems}
\label{sub_pos3}

If $D$ is a dograph
with relations, then
we consider
$\underline{D}$ as the
quiver with relations,
whose \emph{set of
relations consists of
the relations of $D$
and the adjoint
relations} (defined in
\eqref{5.a1a3}).
Suppose we know any
maximal set $\ind
(\underline{D})$ of
nonisomorphic
indecomposable
representations of the
quiver $\underline{D}$
(this means that every
indecomposable
representation of
$\underline{D}$
satisfying the
relations is
isomorphic to exactly
one representation
from $\ind
(\underline{D})$).
Transform
$\ind
(\underline{D})$) as
follows:
\begin{itemize}
  \item
First replace each
representation in
$\ind (\underline{D})$
that is {\it
isomorphic} to a
selfadjoint
representation by one
that is {\it actually}
selfadjoint---i.e.,
has the form
$\underline{\cal A}$,
and denote the set of
these $\underline{\cal
A}$ by
$\ind_0(\underline{D})$.

  \item
Then in each of the
one- or two-element
subsets
$$
\{{\cal M},{\cal
L}\}\subset\ind(\underline{D})
\smallsetminus
\ind_0(\underline{D})\quad
\text{such that }{\cal
M}^{\circ}\simeq {\cal
L},
$$
select one
representation and
denote the set of
selected
representations by
$\ind_1(\underline{D})$.
(If ${\cal M}\sim
{\cal M}^{\circ}$ then
$\{{\cal M},{\cal
L}\}$ consists of one
representation and we
take it.)
\end{itemize}
We obtain a new set
$\ind (\underline{D})$
that we partition into
3 subsets:
\begin{equation}\label{4.8d}
{\ind (\underline{D})}
=
\begin{tabular}{|c|c|}
 \hline &\\[-12pt]
  $\; {\cal M}\;  $&
  ${\cal M}^{\circ}\text{ (if
  ${\cal M}^{\circ}\not
  \simeq{\cal M})$}$\\
  \hline
\multicolumn{2}{|c|}
{$\underline{\cal A}\vphantom{{\hat{N}}}$}\\
 \hline
\end{tabular}\,,\;
\begin{matrix}
 {\cal M}\in
\ind_1(\underline{D}),\\[1pt]
\underline{\cal
A}\in\ind_0(\underline{D}).
\end{matrix}
\end{equation}

For each
representation ${\cal
M}$ of
$\underline{D}$, we
define a
representation ${\cal
M}^+$ of $D$ by
setting ${\cal
M}^+_i:={\cal
M}_i\oplus {\cal
M}_{i^*}^*$ for all
vertices $i$ of $D$
and
\begin{equation}\label{piyf}
{\cal M}^+_{\alpha }:=
\begin{bmatrix} {\cal
M}_{\alpha}&0\\0&{\cal
M}_{\alpha^*}^{*}
 \end{bmatrix},\quad
{\cal M}^+_{\beta}:=
 \begin{bmatrix}
   0&{\cal
M}_{\beta^*}^{*}\\{\cal
M}_{\beta}&0
 \end{bmatrix},\quad
{\cal M}^+_{\gamma}:=
\begin{bmatrix}
   0&{\cal M}_{\gamma}\\
{\cal
M}_{\gamma^*}^{*}&0
 \end{bmatrix}
\end{equation}
for all edges $
\alpha\colon
i\longrightarrow j$,
$\beta\colon i\lin j\
(i\le j)$, and
$\gamma\colon
i\longleftrightarrow
j\ (i\le j).$

The representations
${\cal M}^+$ arise as
follows: each
representation ${\cal
M}$ of $\underline D$,
defines the
selfadjoint
representation ${\cal
M}\oplus {\cal
M}^{\circ}$; the
corresponding
representation of $D$
is ${\cal M}^+$ (and
so $\underline{\cal
M}^+={\cal M}\oplus
{\cal M}^{\circ}$).

For example, if ${\cal
M}$ is the
representation
\eqref{4.5}, then the
selfadjointness of
\begin{equation*}\label{4.9}
{{\cal M}\oplus{\cal
M}^{\circ}:}\qquad
\raisebox{43pt}{\xymatrix@R=20pt@C=50pt{
 {U_1\oplus U_2^{*}}\ar[dd]_{
 \left[\begin{smallmatrix}
   A_1&0\\0&A_2^{*}
 \end{smallmatrix}\right]}
 \ar[ddrr]^(.25){
  \left[\begin{smallmatrix}
   B_1&0\\0&B_2^{*}
 \end{smallmatrix}\right]}&
 &{U_2\oplus U_1^{*}}
  \\
  &\\
 {V_1\oplus V_2^{*}}
 \ar[uurr]^(.75){
  \left[\begin{smallmatrix}
   B_2&0\\0&B_1^{*}
 \end{smallmatrix}\right]}
 \ar@<0.4ex>[rr];[]^{
  \left[\begin{smallmatrix}
   C_2&0\\0&C_1^{*}
 \end{smallmatrix}\right]}
 \ar@<-0.4ex>[rr];[]_{
  \left[\begin{smallmatrix}
   C_1&0\\0&C_2^{*}
 \end{smallmatrix}\right]}
 &&{V_2\oplus V_1^{*}}\ar[uu]_{
  \left[\begin{smallmatrix}
   A_2&0\\0&A_1^{*}
 \end{smallmatrix}\right]}
 }}
\end{equation*}
becomes clear if we
interchange the
summands in each
vector space on the
right, interchanging
respectively the
corresponding strips
in the matrices of
linear mappings:
\begin{equation*}\label{4.10}
\raisebox{43pt}{\xymatrix@R=20pt@C=50pt{
{U_1\oplus
U_2^{*}}\ar[dd]_{
 \left[\begin{smallmatrix}
   A_1&0\\0&A_2^{*}
 \end{smallmatrix}\right]}
 \ar[ddrr]^(.25){
  \left[\begin{smallmatrix}
   0&B_2^{*}\\B_1&0
 \end{smallmatrix}\right]}&
 &{U_1^{*}\oplus U_2}
  \\
&\\
{V_1\oplus V_2^{*}}
 \ar[uurr]^(.75){
  \left[\begin{smallmatrix}
   0&B_1^{*}\\B_2&0
 \end{smallmatrix}\right]}
 \ar@<0.4ex>[rr];[]^{
  \left[\begin{smallmatrix}
   0&C_2\\C_1^{*}&0
 \end{smallmatrix}\right]}
 \ar@<-0.4ex>[rr];[]_{
  \left[\begin{smallmatrix}
   0&C_1\\C_2^{*}&0
 \end{smallmatrix}\right]}
 &&{V_1^{*}\oplus V_2}\ar[uu]_{
  \left[\begin{smallmatrix}
   A_1^{*}&0\\0&A_2
 \end{smallmatrix}\right]}
 }}
\end{equation*}
The corresponding
representation of $D$
is
\[
{\cal M}^+:\qquad
\raisebox{33pt}{
 \xymatrix@R=40pt{
 {U_1\oplus U_2^{*}}\ar[d]_{
 \left[\begin{smallmatrix}
   A_1&0\\0&A_2^{*}
 \end{smallmatrix}\right]}
 \ar@{-}@/^/[d]^{
 \left[\begin{smallmatrix}
   0&B_2^{*}\\B_1&0
 \end{smallmatrix}\right]}
  \\
 {V_1\oplus V_2^{*}}
 &
\save !<-11mm,0cm>
 \ar@(ur,dr)@{<->}^{
 \left[\begin{smallmatrix}
   0&C_1\\C_2^{*}&0
 \end{smallmatrix}\right]}
\restore } }
\]

For every
representation ${\cal
A}$ of ${D}$ and for
every selfadjoint
automorphism
$f=f^{\circ}\colon
\underline{\cal
A}\is\underline{\cal
A}$, we denote by
${\cal A}^f$ the
representation of $D$
obtained from ${\cal
A}$ by replacing
\begin{itemize}
  \item
each form ${\cal
A}_{\beta}$
$(\beta\colon i\lin
j$, $i\le j)$ by
${\cal
A}^f_{\beta}:={\cal
A}_{\beta}f_j$,

  \item
each form ${\cal
A}_{\gamma}$
$(\gamma\colon
i\longleftrightarrow
j$, $i\le j)$ by
${\cal
A}^f_{\gamma}:=f^{-1}_i{\cal
A}_{\gamma}$.
\end{itemize}

The corresponding
selfadjoint
representation
$\underline{{\cal
A}^f}$ of
$\underline{D}$ can be
visualized as the
diagonal of the
rectangle
\begin{equation}\label{ndo}
\begin{split}
 \xymatrix@C=100pt{
 *++[o][F]{v}\ar@{-}[r]
 ^{\underline{\cal A}}
\ar@{--}[rd]^{\underline{{\cal
A}^f}} \ar[d]_{f_v}
 &*++[o][F]{v^*}
 \ar[d]^{f_{v^*}=f_{v}^*} \\
 *++[o][F]{v}
 \ar@{-}[r]^{\underline{\cal A}}
 &*++[o][F]{v^*}
 }
\end{split}
\end{equation}
in which $v$
represents the
vertices of $D$ (thus,
$\underline{\cal
A}\simeq
\underline{{\cal
A}^f}$).

For example, if ${\cal
A}$ is the first
representation in
\eqref{4.2}, then a
selfadjoint
automorphism
\begin{equation*}
\raisebox{57pt}%
{\xymatrix@R=10pt{
 &&{U}\ar@{-->}[ddddddrrr]
 \ar[dddd]_(.30){f_2}
 \ar[ddl]_{A}
 \ar[ddrrr]^(.25){B}&&&
 &{U^{*}}
 \ar[dddd]^{f_2^{*}}
 \\
 {\ \underline{\cal A}:}
 \ar[dddd]_{f}
 \\
 &{V}\ar@{-->}[ddrrrrr]
 \ar[dddd]_{f_1}
 \ar@<0.4ex>[uurrrrr]^(.75){B^{*}}
 \ar@<0.4ex>[rrrr];[]^(.25){C^*}
 \ar@<-0.4ex>[rrrr];[]_(.25){C}
 &&&&{V^{*}}
 \ar[dddd]^(.25){f_1^{*}}
 \ar[uur]_{A^{*}}&
 \\
 \\
 &&{U}\ar[ddl]_{A}
  \ar@<-0.7ex>[ddrrr]_(.25){B}&&&
 &{U^{*}}
 \\
 {\ \underline{\cal A}:}
 \\
 &{V}
 \ar@<0.4ex>[uurrrrr]^(.25){B^{*}}
 \ar@<0.4ex>[rrrr];[]^{C^*}
 \ar@<-0.4ex>[rrrr];[]_{C}
 &&&&{V^{*}}
 \ar@{-->}@<-0.4ex>[uuuullll]
 \ar@{-->}@<0.4ex>[uuuullll]
 \ar[uur]_{A^{*}}&
 }}
 \end{equation*}
defines the
representation
\begin{equation*}
\raisebox{23pt}{\xymatrix@R=4pt{
 &{U}\ar[dd]_{A}
 \ar@{--}@/^/[dd]^{Bf_2}\\
{{\cal A}^{ f}}:&
  \\
 &{V}
\save !<2mm,0cm>
\ar@(ur,dr)%
 @{<-->}^{f_1^{-1}C}
\restore }}
\end{equation*}

Let ${\ind
(\underline{D})}$ be
partitioned as in
\eqref{4.8d}, and let
$\underline{\cal
A}\in{\ind_0
(\underline{D})}$. By
\cite[Lemma
1]{ser_izv}, the set
$R$ of noninvertible
elements of the
endomorphism ring
$\End (\underline{\cal
A})$ is the radical.
Therefore, $\mathbb
T({\cal A}):=\End
(\underline{\cal
A})/R$ is a field or
skew field, on which
we define the
involution
\begin{equation}\label{kyg}
(f+R)^{\circ}:=f^{\circ}+R.
\end{equation}

For each nonzero
$a=a^{\circ}\in
\mathbb T({\cal A})$,
we fix a selfadjoint
automorphism
\begin{equation}\label{lfw}
f_a=f_a^{\circ}\in
a,\quad\text{ and
define ${\cal
A}^a:={\cal A}^{f_a}$}
\end{equation}
(we can take
$f_a:=(f+f^{\circ})/2$
for any $f\in a$). The
set of representations
${\cal A}^a$ is called
the \emph{orbit} of
${\cal A}$. Note that
the corresponding
representations
$\underline{{\cal
A}^a}$ of $\underline
D$ are isomorphic to
$\underline{\cal A}$.
Conversely, if
$\underline{\cal
B}\simeq
\underline{\cal A}$
then ${\cal B}\simeq
{\cal A}^a$ for some
nonzero
$a=a^{\circ}\in
\mathbb T({\cal A})$;
this follows from the
next theorem.

For each Hermitian
form
\[
\varphi(x)=x^{\circ}_1a_1x_1+\dots+
x^{\circ}_ra_rx_r,\qquad
0\ne
a_i=a_i^{\circ}\in
\mathbb T({\cal A}),
\]
we write
\[
{\cal
A}^{\varphi(x)}:=
{\cal
A}^{a_1}\oplus\dots\oplus
{\cal A}^{a_r}.
\]

\begin{theorem}\label{tetete1}
Over a field or skew
field $\mathbb F$ of
characteristic
different from $2$
with involution
$a\mapsto \bar{a}$
$($possibly, the
identity$)$, every
representation of a
dograph $D$ with
relations is
isomorphic to a direct
sum
\begin{equation}\label{iap}
{\cal
M}_1^+\oplus\dots\oplus
{\cal M}_p^+\oplus
{\cal
A}_1^{\varphi_1(x)}\oplus
\dots\oplus {\cal
A}_q^{\varphi_q(x)},
\end{equation}
where
\[
{\cal M}_i\in
\ind_1(\underline{D}),\qquad
\underline{\cal
A}_j\in
\ind_0(\underline{D}),
\]
${\cal A}_j\ne {\cal
A}_{j'}$ if $j\ne j'$,
and each
${\varphi_j(x)}$ is a
Hermitian form over
$\mathbb T({\cal
A}_j)$ with involution
\eqref{kyg}. This sum
is determined by the
original
representation
uniquely up to
permutation of
summands and
replacement of ${\cal
A}_j^{\varphi_j(x)}$
by ${\cal
A}_j^{\psi_j(x)}$,
where ${\psi_j(x)}$ is
a Hermitian form over
$\mathbb T({\cal
A}_j)$ that is
equivalent to
${\varphi_j(x)}$.
\end{theorem}
\begin{proof}
An analogous statement
was proved in
\cite[Theorem
1]{ser_izv} for
selfadjoint
representations of a
linear category with
involution. This
ensures Theorem
\ref{tetete1} since
every dograph $D$ with
relations defines the
following category
$\cal C$ (see
\cite[\S\,2]{ser_izv}):
its objects are the
vertices of
$\underline D$; if $u$
and $v$ are two
vertices of
$\underline D$ then
the set of morphisms
from $u$ to $v$ is the
vector space over the
center of $\mathbb F$
spanned by all
oriented paths from
$u$ to $v$ on
$\underline D$ and
factorized by the
relations on
$\underline D$. An
involution on $\cal C$
is defined in the same
way as the involution
on relations (see
\eqref{5.a1a} and
\eqref{5.a1a3}):
$$
\sum_{i=1}^m
c_i\alpha_{ip_i}\cdots
\alpha_{i2}\alpha_{i1}\
\longmapsto\
\sum_{i=1}^m
\bar{c}_i\alpha_{i1}^*
\alpha_{i2}^*\cdots
\alpha_{ip_i}^*.
\eqno\qedhere
$$
\end{proof}

Theorem \ref{tetete1}
was extended in
\cite{ser_sym} to
symmetric
representations of
algebras with
involution.

For each
representation  $\cal
A$ of $D$, we write
${\cal A}^-:={\cal
A}^{-1}$, where $-1\in
\Aut\underline{\cal
A}$; this means that
the representation
${\cal A}^-$ is
obtained from $\cal A$
by multiplying all the
forms by $-1$:
\begin{equation*}\label{4r}
\raisebox{23pt}{\xymatrix@R=4pt{
 &{U}\ar[dd]_{A}
 \ar@{-}@/^/[dd]^{B}\\
{\cal A}\colon &
  \\
 &{V}
\save !<2mm,0cm>
\ar@(ur,dr)@{<->}^{C}
\restore}}
   \qquad   \qquad   \qquad
\raisebox{23pt}{\xymatrix@R=4pt{
 &{U}\ar[dd]_{A}
 \ar@{-}@/^/[dd]^{-B}\\
{\cal A}^-\colon &
  \\
 &{V}
\save !<2mm,0cm>
\ar@(ur,dr)@{<->}^{-C}
\restore}}
\end{equation*}

Theorem \ref{tetete1}
implies the following
generalization of
Sylvester's Inertia
Theorem.

\begin{theorem}\label{tetete}
Let $\mathbb F$ be
either
\begin{itemize}
  \item[\rm(i)]
an algebraically
closed field of
characteristic
different from $2$
with the identity
involution, or

  \item[\rm(ii)]
an algebraically
closed field with
nonidentity
involution, or

  \item[\rm(iii)]
a real closed field,
or the skew field of
quaternions over a
real closed field.
\end{itemize}

Then every
representation of a
dograph $D$ with
relations over
$\mathbb F$ is
isomorphic to a direct
sum, determined
uniquely up to
permutation of
summands, of
representations of the
following types:
\[
{\cal M}^+,\
  \begin{cases}
   {\cal
A}& \text{if ${\cal
A}^{-}\simeq{\cal
A}$}, \\
{\cal A},\ {\cal
A}^{-} & \text{if
${\cal
A}^{-}\not\simeq{\cal
A}$},
  \end{cases}\quad(
\text{where } {\cal
M}\in
\ind_1(\underline{P})
,\ \underline{\cal
A}\in\ind_0(\underline{P}))),
\]
or, respectively to
the cases
{\rm(i)--(iv)},
\begin{itemize}
  \item[\rm(i)]
${\cal M}^+$, ${\cal
A}$,

  \item[\rm(ii)]
${\cal M}^+$, ${\cal
A}$, ${\cal A}^{-}$,

  \item[\rm(iii)]
$ {\cal M}^+,
  \begin{cases}
    \ \ {\cal A}, &
\parbox[t]{9cm}{if $\mathbb
T({\cal A})$ is an
algebraically closed
field with the
identity involution or
a skew field of
quaternions with
involution different
from quaternionic
conjugation, and} \\
    {\cal A},{\cal A}^{-}, &
    \text{otherwise}.
  \end{cases}
$
\end{itemize}
\end{theorem}

\begin{proof}
Theorem \ref{tetete1}
reduces the
classification of
representations of any
dograph $D$ to the
classification of
Hermitian forms over
the fields or skew
fields $\mathbb
T({\cal A})$,
$\underline{\cal
A}\in\ind_0(\underline{D})$,
assuming known
$\ind_1(\underline{D})$,
$\ind_0(\underline{D})$,
and the orbit of the
representations ${\cal
A}$ for each
$\underline{\cal
A}\in\ind_0(\underline{D})$.

If $\mathbb F$ is
finite dimensional
over its center
$C(\mathbb F)$, then
$\mathbb T({\cal A})$
is also finite
dimensional over
$C(\mathbb F)$ under
the natural imbedding
of $C(\mathbb F)$ into
the center of $\mathbb
T({\cal A})$, and the
involution on $\mathbb
T({\cal A})$ extends
the involution on
$C(\mathbb F)$.
\medskip

(i) If $\mathbb F$ is
an algebraically
closed field of
characteristic
different from 2 with
the identity
involution, then it
has no finite
extensions. Hence,
$\mathbb T(\cal
A)=\mathbb F$ for each
$\underline{\cal
A}\in\ind_0(\underline{D})$,
and so each Hermitian
form
\begin{equation}\label{ndir}
a_1x_1^2+\dots+a_rx_r^2,\qquad
0\ne a_i\in\mathbb F,
\end{equation}
is equivalent to
$x_1^2+\dots+x_r^2$.
We can replace all
${\cal
A}_j^{\varphi_j(x)}$
in \eqref{iap} by
${\cal
A}_j\oplus\dots\oplus
{\cal A}_j$. In view
of Theorem
\ref{tetete1}, the
obtained direct sum is
determined by the
original
representation
uniquely up to
permutation of
summands.
\medskip

(ii) Let $\mathbb F$
be an algebraically
closed field with
nonidentity
involution. Its
characteristic is $0$
by Lemma \ref{l00}(b),
$\mathbb T(\cal
A)=\mathbb F$ for each
$\underline{\cal
A}\in\ind_0(\underline{D})$,
and the involution
$a\mapsto a^{\circ}$
on $\mathbb T(\cal A)$
coincides with the
involution $a\mapsto
\bar a$ on $\mathbb
F$. Due to the law of
inertia \cite{bour2},
each Hermitian form
\[
a_1\bar
x_1x_1+\dots+a_r\bar
x_rx_r,\qquad 0\ne
a_i=\bar a_i\in\mathbb
F,
\]
is equivalent to
exactly one form
\[
\bar x_1x_1+\dots
+\bar x_lx_l -\bar
x_{l+1}x_{l+1} -\dots-
\bar x_rx_r.
\]
Therefore, we can
replace each ${\cal
A}_j^{\varphi_j(x)}$
in \eqref{iap} by
exactly one direct sum
of the form
\begin{equation}\label{bdp}
{\cal
A}_j\oplus\dots\oplus
{\cal A}_j\oplus {\cal
A}_j^-\oplus\dots\oplus
{\cal A}_j^-.
\end{equation}
\medskip

(iii) Let $\mathbb F$
be a real closed field
$\mathbb P$ or the
skew field $\mathbb H$
of quaternions over a
real closed field
$\mathbb P$. By Lemma
\ref{l00}(a),
$\charact (\mathbb
F)=0$. The center of
$\mathbb F$ is
$\mathbb P$. Hence,
$\mathbb T(\cal A)$
for each
$\underline{\cal
A}\in\ind_0(\underline{D})$
is a finite extension%
\footnote{Formulating
Theorem 2 in
\cite{ser_izv}, I
erroneously thought
that all $\mathbb
T(\cal A)=\mathbb H$
if $\mathbb F=\mathbb
H$. To correct it,
remove ``or the
algebra of quaternions
\dots'' in a) and b)
and add ``or the
algebra of quaternions
over a maximal ordered
field'' in c). The
paper \cite{ser1} is
based on the incorrect
Theorem 2 in
\cite{ser_izv} and so
the signs $\pm$ of the
sesquilinear forms in
the indecomposable
direct summands in
\cite[Theorems
1--4]{ser1} are
incorrect. Correct
canonical forms are
given for pairs of
symmetric/skew-symmetric
matrices in
\cite{rod_pair_nonst,
rod_pair_stand,
rod_pair_herm}, for
selfadjoint operators
in \cite{kar}, and for
isometries in Theorem
\ref{theor}.}\label{page}
of $\mathbb P$. By the
Frobenius theorem
\cite{bour2}, $\mathbb
T(\cal A)$ is either
$\mathbb P$, or its
algebraic closure
$\mathbb P+\mathbb
Pi$, or $\mathbb H$.

If $\mathbb T(\cal A)$
is either $\mathbb P$,
or $\mathbb P+\mathbb
Pi$ with nonidentity
involution, or
$\mathbb H$ with
quaternionic
conjugation
\eqref{ne}, then by
the law of inertia
\cite{bour2} each
Hermitian form
\begin{equation}\label{ksre}
\varphi(x)=x_1^{\circ}a_1x_1+\dots
+x_r^{\circ}a_rx_r,\qquad
0\ne
a_i=a_i^{\circ}\in\mathbb
T(\cal A),
\end{equation}
is equivalent to
exactly one form
\[
x_1^{\circ}x_1+\dots
+x_l^{\circ}x_l
-x_{l+1}^{\circ}x_{l+1}
-\dots-
x_r^{\circ}x_r.
\]
Therefore, we can
replace each ${\cal
A}_j^{\varphi_j(x)}$
in \eqref{iap} by
\eqref{bdp}.

If $\mathbb T({\cal
A}) =\mathbb P+\mathbb
Pi$ with the identity
involution, then every
Hermitian form over it
is equivalent to
$x_1^2+\dots+x_r^2$.
We can replace each
${\cal
A}_j^{\varphi_j(x)}$
in \eqref{iap} by
${\cal
A}_j\oplus\dots\oplus
{\cal A}_j$.

If $\mathbb T({\cal
A}) = \mathbb H$ with
quaternionic
semiconjugation
\eqref{nen}, then
every Hermitian form
\eqref{ksre} is
equivalent to $
x_1^{\circ}x_1+\dots
+x_r^{\circ}x_r $
since each $a_i$ is
represented in the
form $a_i =
b_i^{\circ}b_i$, where
$b_i:=\sqrt{a_i}$ is
taken in the field
$\mathbb P+\mathbb Pj$
if $a_i\in\mathbb P$,
or in the field
$\mathbb P(a_i)$ if
$a_i\notin\mathbb P$;
these fields are
algebraically closed
and the involution
\eqref{nen} acts
identically on them.
Therefore, we can
replace each ${\cal
A}_j^{\varphi_j(x)}$
in \eqref{iap} by
${\cal
A}_j\oplus\dots\oplus
{\cal A}_j$.
\end{proof}

\begin{example}
The problem of
classifying Hermitian
forms over $\mathbb F$
is given by the
dograph
\[
\xymatrix{ D\,: &{1}
 \ar@(ur,dr)@{-}^{\alpha}
 }\qquad\qquad
 \alpha^*=\alpha.
\]
Its quiver is
\[
\underline{
D}:\quad\xymatrix{
 {1}
  \ar@/^/@{->}[rr]^{\alpha}
 \ar@/_/@{->}[rr]_{\alpha^*} &&{1^*}
 }\qquad\qquad
 \alpha^*=\alpha.
\]
Each representation
\[
{\cal
M}:\quad\xymatrix{
 {U}
  \ar@/^/@{->}[rr]^A
 \ar@/_/@{->}[rr]_B
 &&{V}
 }\qquad\qquad
 A=B
\]
of $\underline{D}$ is
given by a linear
mapping $A\colon U\to
V$, which is a direct
sum of mappings of the
types $\mathbb
F\stackrel{1}{\to}\mathbb
F$, $0\to\mathbb F$,
and $\mathbb F\to 0$
(because each matrix
reduces by equivalence
transformations to
$I_r\oplus 0$, which
is a direct sum of
matrices of the types
$I_1$, $0_{10}$, and
$0_{01}$). Thus, the
set \eqref{4.8d} is
\[
{\ind (\underline{D})}
=
\begin{tabular}{|c|c|}
 \hline &\\[-12pt]
  $\; 0\rightrightarrows
  {\mathbb F}\;  $&
  $\; {\mathbb F}
  \rightrightarrows
  0\;  $\\
  \hline
\multicolumn{2}{|c|}
{$\; {\mathbb F}
  \rightrightarrows
  {\mathbb F}\;  $}\\
 \hline
\end{tabular}\,.
\]
Theorem \eqref{tetete}
ensures that every
representation of $D$
is isomorphic to a
direct sum of
representations of the
types
$\xymatrix{{{\mathbb
F}}
 \ar@(ur,dr)@{-}^{0}
 }$ (which is $[0\rightrightarrows
  {\mathbb F}]^+$)
and
$\xymatrix{{{\mathbb
F}}
 \ar@(ur,dr)@{-}^{1}
 }$ if $\mathbb F$ is an
algebraically closed
field of
characteristic not $2$
with the identity
involution, and of
representations of the
types
$\xymatrix{{{\mathbb
F}}
 \ar@(ur,dr)@{-}^{0}
 }$,
$\xymatrix{{{\mathbb
F}}
 \ar@(ur,dr)@{-}^{1}
 }$,
and
$\xymatrix{{{\mathbb
F}}
 \ar@(ur,dr)@{-}^{-1}
 }$ if $\mathbb F$ is an
algebraically closed
field with nonidentity
involution.
\end{example}

\begin{corollary}
Each system of linear
mappings and
bilinear/sesquilinear
forms on vector spaces
over $\mathbb R$,
$\mathbb C$, or
$\mathbb H$ decomposes
into a direct sum of
indecomposable systems
uniquely up to
isomorphisms of
summands.
\end{corollary}

By \cite{ser_disch},
the set of dimensions
\eqref{bst} of
indecomposable
representations of a
dograph does not
depend on the
orientation of its
edges, and so by Kac's
theorem \cite{kac} it
coincides with the set
of positive roots of
the dograph. An
analogous description
of the set of
dimensions of
indecomposable
Euclidean or unitary
representations of a
quiver (i.e., each
vertex is assigned by
a Euclidean or unitary
space) is given in
\cite{ser_unit}.

\section{Proof of
Theorem \ref{Theorem
5}} \label{s_pro}

Each pair $(A,B)$
consisting of a
nondegenerate
$\varepsilon$-Hermitian
form $B$ and an
isometric operator $A$
on a vector space $V$
over a field or skew
field of
characteristic
different from $2$
determines the
representation
\[
\xymatrix{{\cal A} :&&
 {V\,}
\save !<-2mm,0cm>
\ar@(ul,dl)@{->}_{A}
\restore
 \ar@<-0.4ex>@(u,r)
@{-}^{B}
\ar@<-0.4ex>@(r,d)
@{<->}^{B^{-1}}
  }
\]
of the dograph $D$
defined in
\eqref{2.8akda}; if
$B$ is given by a
matrix $B_e$ in some
basis of $V$ then
$B^{-1}$ is given by
$B_e^{-1}$ in the
*dual basis of $V^*$.

The quiver with
involution of the
dograph $D$ is
\begin{equation*}
\underline
D:\quad\xymatrix{
 {1}\ar@(ul,dl)@{->}_{\alpha}
 \ar@/^2pc/@{->}[rr]^{\beta}
  \ar@/^/@{->}[rr]^{\beta^*}
 \ar@/^/@{->}[rr];[]^{\gamma}
 \ar@/^2pc/@{->}
 [rr];[]^{\gamma^*} &&{1^*}
\save !<1.5mm,0cm>
\ar@(ur,dr)@{->}^{\alpha^*}
\restore
 }\qquad
 {\begin{matrix}
 \beta=\varepsilon \beta^*
 =\alpha^*\beta\alpha,\\
 \gamma\beta=1,\quad
 \beta\gamma=1,\\
 \gamma^*\beta^*=1,\quad
 \beta^*\gamma^*=1.
 \end{matrix}}
\end{equation*}

We will prove Theorem
\ref{Theorem 5} using
Theorem \ref{tetete1};
to do this, we first
identify in Lemma
\ref{lenhi} the sets
$\ind_1(\underline{D})$
and
$\ind_0(\underline{D})$,
and the orbit of $\cal
A$ for each
$\underline{\cal A}\in
\ind_0(\underline{D})$.

The arrow $\gamma$ was
appended in $D$ with
the only purpose: each
form assigned to
$\beta $ must be
nondegenerate. So we
will omit $\gamma$ and
$\gamma^*$ and
represent $D$ and
$\underline D$ as
follows:
\begin{equation}\label{jsos+}
{D}:\quad\xymatrix{
 {1} \ar@(ul,dl)@{->}_{\alpha}
 \ar@(ur,dr)@{-}^{\beta}}
\qquad
 \begin{matrix}
 \beta=\varepsilon \beta^*
 =\alpha^*\beta\alpha,\\
\beta \text{ is
nonsingular,}
\end{matrix}
\end{equation}
\begin{equation}\label{jso+}
\underline{D}:\quad\xymatrix{
 {1} \ar@(ul,dl)@{->}_{\alpha}
  \ar@/^/@{->}[rr]^{\beta}
 \ar@/_/@{->}[rr]_{\beta^*} &&{1^*}
\save !<1.5mm,0cm>
\ar@(ur,dr)@{->}^{\alpha^*}
\restore
 }\qquad
 \begin{matrix}
 \beta=\varepsilon \beta^*
 =\alpha^*\beta\alpha,\\
\beta \text{ is
bijective.}
\end{matrix}
\end{equation}

Every representation
of $D$ or $\underline
D$ over $\mathbb F$ is
isomorphic to a
representation in
which all vector
spaces are $\mathbb
F\oplus\dots\oplus\mathbb
F$. We will consider
only such
representations of $D$
and $\underline D$,
they can be given by
matrix pairs $(A,B)$:
\begin{equation*}\label{ren7a+}
{\cal A}:\quad
\xymatrix{
 *{\ci}
\ar@(ul,dl)_{A}
 \ar@(ur,dr)@{-}^{
B}} \qquad
 \begin{matrix}
B=\varepsilon B^*=A^*BA,\\
B \text{ is
nonsingular,}
\end{matrix}
 \end{equation*}
and, respectively, by
matrix triples ${\cal
M}=(A,B,C)$:
\begin{equation}\label{ser14sd}
{\cal
M}:\quad\xymatrix{
 *{\ci} \ar@(ul,dl)@{->}_{A}
  \ar@/^/@{->}[rr]^{B}
 \ar@/_/@{->}[rr]_{\varepsilon B}
 &&*{\ci}
\ar@(ur,dr)@{->}^{C}
 }\qquad
 \begin{matrix}
B= CBA,\\
B \text{ is
nonsingular}
\end{matrix}
    \end{equation}
(we omit the spaces
$\mathbb
F\oplus\dots\oplus\mathbb
F$ since they are
completely determined
by the sizes of the
matrices).

The adjoint
representation
\begin{equation*}\label{ren2+}
{\cal
M}^{\circ}:\quad\xymatrix{
 *{\ci} \ar@(ul,dl)@{->}_{C^*}
  \ar@/^/@{->}[rr]^{\varepsilon
  B^*}
 \ar@/_/@{->}[rr]_{B^*}
 &&*{\ci}
\ar@(ur,dr)@{->}^{A^*}
 }
 \end{equation*}
is given by the matrix
triple
\begin{equation}\label{ldbye}
(A,B,C)^{\circ} =
(C^*,\varepsilon
B^*,A^*).
\end{equation}

A morphism of
representations
\begin{equation*}\label{ser14d}
\xymatrix@R=12pt{{\
{\cal M}:}\ar[dd]_g&&
 *{\ci}\ar[dd]_{G_1}
 \ar@(ul,dl)@{->}_{A}
  \ar@/^/@{->}[rr]^{B}
 \ar@/_/@{->}[rr]_{\varepsilon B}
 &&*{\ci}\ar[dd]^{G_2}
\ar@(ur,dr)@{->}^{C}
\\ \\
 {\ {\cal M}':}&&*{\ci}
 \ar@(ul,dl)@{->}_{A'}
  \ar@/^/@{->}[rr]^{B'}
 \ar@/_/@{->}[rr]_{\varepsilon B'}
 &&*{\ci}
\ar@(ur,dr)@{->}^{C'}
 }
\end{equation*}
is given by the matrix
pair
\[
g=[G_1,G_2]\colon\;
{\cal M}\to {\cal M}'
\]
(for morphisms we use
square brackets)
satisfying
\begin{equation}\label{msp}
G_1A=A'G_1,\quad
G_2B=B'G_1,\quad
G_2C=C'G_2,
\end{equation}
and the adjoint
morphism is given by
the matrix pair
\begin{equation}\label{ami}
g^{\circ}=[G_2^*,G_1^*]
\colon\; {\cal
M}^{\prime\circ}\to
{\cal M}^{\circ}.
\end{equation}

\begin{lemma}\label{lenhi}
Let $\mathbb F$ be a
field or skew field of
characteristic
different from $2$.
Let ${\cal O}_{\mathbb
F}$ be a maximal set
of nonsingular
indecomposable
canonical matrices
over\/ $\mathbb F$ for
similarity $($see page
{\rm\pageref{papage}}$)$.
Let $D$ be the dograph
\eqref{jsos+}. Then
the following
statements hold:

\begin{itemize}
  \item[{\rm(a)}] The set
$\ind(\underline{{D}})$
can be taken to be all
representations
$(\Phi, I, \Phi^{-1})$
with $\Phi\in{\cal
O}_{\mathbb F}$.

\item[\rm(b)] The set
$\ind_1(\underline{{D}})$
can be taken to be all
representations
$${\cal M}_{\Phi}:=
(\Phi,I,\Phi^{-1}),$$
in which $\Phi\in{\cal
O}_{\mathbb F}$ is
such that
$\Phi_{(\varepsilon)}$
$($defined in
\eqref{mshd}$)$ does
not exist, and
\begin{equation}\label{4.adg}
\parbox{15em}
{$\Phi$ is determined
up to replacement by
$\Psi\in{\cal
O}_{\mathbb F}$ that
is similar to
$\Phi^{-*}$.}
\end{equation}
The corresponding
representation
\eqref{piyf} of ${D}$
has the form
\begin{equation}\label{ren5mi}
{\cal M}_{\Phi}^+:
\quad\xymatrix{
 *{\ci}
\ar@(ul,dl)_{
\begin{bmatrix}
 \Phi&0\\0&\Phi^{-*}
\end{bmatrix}
}
 \ar@(ur,dr)@{-}^{
\begin{bmatrix}
0&\varepsilon I\\
 I&0
\end{bmatrix}.
}}
\end{equation}

  \item[\rm(c)]
The set
$\ind_0(\underline{{D}})$
can be taken to be all
representations
\begin{equation}\label{nsp}
\underline{\cal
A}_{\Phi} :=
(\Phi,\Phi_{(\varepsilon)},\Phi^*),
\end{equation}
in which $\Phi\in{\cal
O}_{\mathbb F}$ is
such that
$\Phi_{(\varepsilon)}$
exists. The
corresponding
representations of
${D}$ have the form
\begin{equation*}\label{5mi}
{\cal A}_{\Phi}: \ \
\xymatrix{
 *{\ci}
\ar@(ul,dl)_{\Phi}
 \ar@(ur,dr)@{-}^{
\Phi_{(\varepsilon)}}},
\qquad {\cal
A}_{\Phi}^-: \ \
\xymatrix{
 *{\ci}
\ar@(ul,dl)_{\Phi}
 \ar@(ur,dr)@{-}^{
-\Phi_{(\varepsilon)}}}.
\end{equation*}

\item[\rm(d)] Let
$\mathbb F$ be a field
and let
$\underline{\cal
A}_{\Phi} :=
(\Phi,\Phi_{(\varepsilon)},\Phi^*)
\in\ind_0(\underline{{D}})$.

\begin{itemize}
  \item[\rm(i)]
The ring
$\End(\underline{\cal
A}_{\Phi})$ of
endomorphisms of
$\underline{\cal
A}_{\Phi}$ consists of
the matrix pairs
\begin{equation}\label{ldy}
[f(\Phi),f(\Phi^{-*})],\qquad
f(x)\in\mathbb F[x],
\end{equation}
and the involution on
$\End(\underline{\cal
A}_{\Phi})$ is
\[
[f(\Phi),f(\Phi^{-*})]^{\circ}=
[\bar
f(\Phi^{-1}),\bar
f(\Phi^*)].
\]

  \item[\rm(ii)]
$\mathbb T({\cal
A}_{\Phi})$ can be
identified with the
field
\begin{equation}\label{ksy}
{\mathbb
F}(\kappa)={\mathbb
F}[x]/p_{\Phi}(x){\mathbb
F}[x],\qquad
\kappa:=x+p_{\Phi}(x){\mathbb
F}[x],
\end{equation}
$(p_{\Phi}(x)$ is
defined in
\eqref{ser24}$)$ with
involution
\begin{equation}\label{kxu}
f(\kappa)^{\circ}=
\bar f(\kappa^{-1}).
\end{equation}
Each element of
$\mathbb T({\cal
A}_{\Phi})$ on which
this involution acts
identically is
uniquely represented
in the form
$q(\kappa)$ for some
nonzero function
\eqref{ser13}. The
representations
\begin{equation}\label{wmp}
{\cal
A}_{\Phi}^{q(\kappa)}:
\quad \ \xymatrix{
 *{\ci}
\ar@(ul,dl)_{\Phi}
 \ar@(ur,dr)@{-}^{
\Phi_{(\varepsilon)}q(\Phi)}}
\end{equation}
$($see \eqref{lfw}$)$
constitute the orbit
of ${\cal A}_{\Phi}$.
\end{itemize}
\end{itemize}
\end{lemma}

\begin{proof}
(a) Every
representation of the
quiver
$\underline{{D}}$ is
isomorphic to one of
the form $(A, I,
A^{-1})$. By
\eqref{msp}, $(A, I,
A^{-1})\simeq (B, I,
B^{-1})$ if and only
if the matrices $A$
and $B$ are similar.

(b)\&(c) Let
$\Phi,\Psi\in{\cal
O}_{\mathbb F}$. In
view of \eqref{ldbye},
\begin{equation}\label{mdpiug}
(\Psi,I,\Psi^{-1})\simeq
(\Phi,I,\Phi^{-1})^{\circ}
=
(\Phi^{-*},\varepsilon
I,\Phi^{*})
\end{equation}
if and only if $\Psi$
is similar to
$\Phi^{-*}$.

Suppose
$(\Phi,I,\Phi^{-1})$
is isomorphic to a
selfadjoint
representation:
\begin{equation}\label{bgd}
[G_1,G_2]\colon\;(\Phi,I,\Phi^{-1})
\is (C,D,C^*),\qquad
D=\varepsilon D^*.
\end{equation}
Define a selfadjoint
representation
$(A,B,A^*)$ by the
congruence
\begin{equation}\label{wfw}
[G_1^{-1},G_1^*]\colon\;
(C,D,C^*) \is
(A,B,A^*),\qquad
B=\varepsilon B^*.
\end{equation}
The composition of
\eqref{bgd} and
\eqref{wfw} is the
isomorphism
\[
\xymatrix@R=15pt{
 *{\ci}\ar[dd]_{I}
 \ar@(ul,dl)@{->}_{\Phi}
  \ar@/^/@{->}[rr]^{I}
 \ar@/_/@{->}[rr]_{\varepsilon I}
 &&*{\ci}\ar[dd]^{G:=G^*_1G_2}
\ar@(ur,dr)@{->}^{\Phi^{-1}}
\\ \\
 *{\ci}
 \ar@(ul,dl)@{->}_{A}
  \ar@/^/@{->}[rr]^{B}
 \ar@/_/@{->}[rr]_{\varepsilon B}
 &&*{\ci}
\ar@(ur,dr)@{->}^{A^*}
 }
\]
By \eqref{msp}, $ A =
\Phi,$ $B = G,$ and
$A^*G = G\Phi^{-1}$;
hence $B=\varepsilon
B^*=\Phi^*B\Phi$. We
can replace $B$ by
$\Phi_{(\varepsilon)}$
and obtain
\[
[I,\Phi_{(\varepsilon)}]\colon
(\Phi,I,\Phi^{-1})\is
(\Phi,\Phi_{(\varepsilon)},\Phi^*).
\]

This means that if
$(\Phi,I,\Phi^{-1})\in
\ind(\underline{{D}})$
is isomorphic to a
selfadjoint
representation then it
is isomorphic to
\eqref{nsp}. Hence the
representations
\eqref{nsp} form
$\ind_0(\underline{{D}})$.
Due to \eqref{mdpiug},
we can identify
isomorphic
representations in the
set of remaining
representations
$(\Phi,I,\Phi^{-1})\in
\ind(\underline{{D}})$
by imposing the
condition
\eqref{4.adg}, and
obtain
$\ind_1(\underline{{D}})$.

(d) Let $\mathbb F$ be
a field. It is known
that if $\Phi$ is a
square matrix over
$\mathbb F$ being
indecomposable for
similarity, then each
matrix over $\mathbb
F$ commuting with
$\Phi$ is a polynomial
in $\Phi$. Let us
recall the proof. We
can assume that $\Phi$
is an $n\times n$
Frobenius block
\eqref{3.lfo}. Then
the vectors
\begin{equation}\label{gwz}
e:=(1,0,\dots,0)^T,\
\Phi e,\ \dots,\
\Phi^{n-1} e
\end{equation}
 form a
basis of $\mathbb
F^n$. Let $S\in\mathbb
F^{n\times n}$ commute
with $\Phi$ and let
\[S e=a_0e+a_1\Phi
e+\dots+
a_{n-1}\Phi^{n-1}e,\qquad
a_0,\dots,
a_{n-1}\in\mathbb F.\]
Define
\[f(x):=a_0+a_1x
+\dots+
a_{n-1}x^{n-1}\in\mathbb
F[x].\] Then
$Se=f(\Phi)e,$
\[S\Phi e=\Phi Se=\Phi f(\Phi)e=
f(\Phi)\Phi e, \
\dots, \ S\Phi^{n-1}
e=f(\Phi)\Phi^{n-1}e.
\]
Since \eqref{gwz} is a
basis, $S=f(\Phi)$.

(i) Let $g=[G_1,
G_2]\in
\End(\underline{\cal
A}_{\Phi}),$ where
$\underline{\cal
A}_{\Phi} =
(\Phi,\Phi_{(\varepsilon)},
\Phi^*)\in\ind_0(\underline{{D}})$.
Then by \eqref{msp}
\begin{equation}\label{mdtc}
G_1\Phi=\Phi G_1,\quad
G_2\Phi_{(\varepsilon)}
=\Phi_{(\varepsilon)}
G_1,\quad
G_2\Phi^*=\Phi^* G_2.
\end{equation}
By \eqref{mshd},
\begin{equation}\label{kst}
\Phi^{-*}=
\Phi_{(\varepsilon)}
\Phi
\Phi_{(\varepsilon)}^{-1}.
\end{equation}
Since $G_1$ commutes
with $\Phi$, we have
\begin{equation*}\label{lrsh}
G_1 = f(\Phi)\ \
(f(x)\in \mathbb
F[x]),\qquad G_2 =
\Phi_{(\varepsilon)}
f(\Phi)
\Phi_{(\varepsilon)}^{-1}=
f(\Phi^{-*}).
\end{equation*}

Consequently, the ring
$\End(\underline{\cal
A}_{\Phi})$ of
endomorphisms of
$\underline{\cal
A}_{\Phi}$ consists of
the matrix pairs
\eqref{ldy}, and the
involution
\eqref{kdtc} has the
form:
\[
[f(\Phi),f(\Phi^{-*})]^{\circ}=
[f(\Phi^{-*})^*,f(\Phi)^*]=
[\bar
f(\Phi^{-1}),\bar
f(\Phi^*)].
\]

(ii) Since
$\Phi_{(\varepsilon)}$
is fixed and $G_2 =
\Phi_{(\varepsilon)}
f(\Phi)
\Phi_{(\varepsilon)}^{-1}$,
each endomorphism
$[f(\Phi),f(\Phi^{-*})]$
is completely
determined by
$f(\Phi)$, and so
$\End(\underline{\cal
A}_{\Phi})$ can be
identified with the
ring
$$ \mathbb
F[\Phi]=\{f(\Phi)\,|\,f\in\mathbb
F[x]\}\quad \text{with
involution
$f(\Phi)\mapsto \bar
f(\Phi^{-1})$,}$$
which is isomorphic to
$\mathbb
F[x]/p_{\Phi}(x)^s\mathbb
F[x]$, where
$p_{\Phi}(x)^s$ is the
characteristic
polynomial
\eqref{ser24} of
$\Phi$. The radical of
$\End(\underline{\cal
A}_{\Phi})$ is
generated by
$p_{\Phi}(\Phi)$,
hence $\mathbb T({\cal
A}_{\Phi})$ is
naturally isomorphic
to the field
\eqref{ksy} with
involution
$f(\kappa)^{\circ}=
\bar f(\kappa^{-1})$.

According to Lemma
\ref{LEMMA 7}, each
element of this field
on which the
involution acts
identically is
uniquely representable
in the form
$q(\kappa)$ for some
nonzero function
\eqref{ser13}.

The pair
$[q(\Phi),\Phi_{(\varepsilon)}
q(\Phi)
\Phi_{(\varepsilon)}^{-1}]$
is an endomorphism of
$\underline{\cal
A}_{\Phi}$ due to
\eqref{mdtc}. This
endomorphism is
selfadjoint since the
function \eqref{ser13}
fulfils
$q(x^{-1})=\bar q(x)$,
and so by \eqref{kst}
\[
\Phi_{(\varepsilon)}
q(\Phi)
\Phi_{(\varepsilon)}^{-1}
=q(\Phi^{-*})=\bar
q(\Phi^*)=
q(\Phi)^*.\]

Since distinct
functions $q(x)$ give
distinct $q(\kappa)$
and \[q(\Phi)\in
q(\kappa)=q(\Phi)
+p_{\Phi}(\Phi){\mathbb
F}[\Phi],\] we can
take in \eqref{lfw} \[
f_{q(\kappa)}:=
[q(\Phi),q(\Phi)^*]
\in
\End(\underline{\cal
A}_{\Phi}).\] By
\eqref{ndo}, the
corresponding
representations ${\cal
A}_{\Phi}^{q(\kappa)}
= {\cal
A}_{\Phi}^{f_{q(\kappa)}}
$ have the form
\eqref{wmp} and
constitute the orbit
of ${\cal A}_{\Phi}$.
\end{proof}

\begin{proof}[Proof of Theorem
\ref{Theorem 5}] Each
pair $(A,B)$
consisting of a
nondegenerate
$\varepsilon$-Hermitian
form $B$ and an
isometric operator $A$
gives a representation
of the dograph
\eqref{jsos+}. By
Theorem \ref{tetete1},
each representation of
\eqref{jsos+} over a
field $\mathbb F$ of
characteristic
different from $2$ is
isomorphic to a direct
sum of representations
of the form ${\cal
M}^+$ and ${\cal
A}^a$, where ${\cal
M}\in
\ind_1(\underline{D})$,
$\underline{\cal A}\in
\ind_0(\underline{D})$,
and $0\ne
a=a^{\circ}\in\mathbb
T({\cal A})$. This
direct sum is
determined uniquely up
to permutation of
summands and
replacement of the
whole group of
summands $ {\cal
A}^{a_1}
\oplus\dots\oplus
{\cal A}^{a_s} $ with
the same $\cal A$ by $
{\cal A}^{b_1}
\oplus\dots\oplus
{\cal A}^{b_s} $ such
that the Hermitian
forms $
a_1x_1^{\circ}x_1+\dots+
a_sx_s^{\circ}x_s$ and
$
b_1x_1^{\circ}x_1+\dots+
b_sx_s^{\circ}x_s$ are
equivalent over the
field $\mathbb T({\cal
A})$ (see
\eqref{ksy}).

This proves Theorem
\ref{Theorem 5} since
we can use the sets
$\ind_1(\underline{D})$
and
$\ind_0(\underline{D})$
from Lemma
\ref{lenhi}, the field
$\mathbb T({\cal A})$
is determined in
\eqref{ksy}, and the
representations ${\cal
M}^+$ and ${\cal A}^a$
have the form
\eqref{ren5mi} and
\eqref{wmp}.
\end{proof}

\section{Proof of Theorem
\ref{theor}}
\label{secmat}

Theorem \ref{theor}
gives canonical
matrices of
representations of the
dograph \eqref{jsos+}
over algebraically or
real closed fields and
over skew fields of
quaternions. We will
prove it basing on the
next lemma, in which
we concretize Lemma
\ref{lenhi}: we give
the sets ${\cal
O}_{\mathbb F}$,
establish when $\Psi$
is similar to
$\Phi^{-*}$ (see
\eqref{4.adg}) and
when
$\Phi_{(\varepsilon)}$
exists for
$\Phi\in{\cal
O}_{\mathbb F}$,
construct the matrices
$\Phi_{(\varepsilon)}$
simpler than in Lemma
\ref{lsdy1}, and find
the field $\mathbb
T({\cal A})$ for each
$\underline{\cal A}\in
\ind_0(\underline{D})$.

Recall the $n$-by-$n$
matrices defined in
\eqref{ase} and
\eqref{ases}:
\[
 \Lambda_n=\begin{bmatrix}
1&2&\cdots&2
\\&1&\ddots&\vdots
\\&&\ddots&2
\\0&&&1
\end{bmatrix},
\qquad
F_n=\begin{bmatrix}
0&&&\ddd&
\\&&1&
\\&-1&\\
1&&&0
\end{bmatrix}.
\]

\begin{lemma}
\label{lsdy}
\begin{itemize}
  \item[\rm(a)]
  Let $\mathbb F$ be an
algebraically closed
field of
characteristic
different from $2$
with the identity
involution, and let
$\varepsilon=\pm 1$.
One can take ${\cal
O}_{\mathbb F}$ to be
all nonsingular Jordan
blocks.  For nonzero
$\lambda,\mu\in\mathbb
F$,
\begin{equation*}\label{slust}
J_n(\lambda)\text{ is
similar to
}J_n(\mu)^{-T}
\quad\Longleftrightarrow\quad
\lambda={\mu}^{-1},
\end{equation*}
\begin{equation*}\label{nslsi}
J_n(\lambda)_{(\varepsilon)}\
\text{ exists
}\quad\Longleftrightarrow\quad
\lambda=\pm 1\text{
and }
\varepsilon=(-1)^{n+1}.
\end{equation*}
If it exists then
$J_n(\lambda)$ is
similar to
\begin{equation}\label{ksidq}
\Psi=\lambda\Lambda_n,\qquad
\text{with \ }
\Psi_{(\varepsilon)}=F_n.
\end{equation}

  \item[\rm(b)]
Let $\mathbb F$ be an
algebraically closed
field with nonidentity
involution. One can
take ${\cal
O}_{\mathbb F}$ to be
all nonsingular Jordan
blocks. For nonzero
$\lambda,\mu\in\mathbb
F$,
\begin{equation*}\label{lust2}
J_n(\lambda)\text{ is
similar to
}J_n(\mu)^{-*}
\quad\Longleftrightarrow\quad
\lambda=\bar{\mu}^{-1},
\end{equation*}
\begin{equation*}\label{nlsi}
J_n(\lambda)_{(1)}\
\text{ exists
}\quad\Longleftrightarrow\quad
|\lambda|=1\quad
(\text{see
\eqref{1kk}}).
\end{equation*}
If it exists then
$J_n(\lambda)$ is
similar to
\begin{equation}\label{ksim}
\Psi=\lambda\Lambda_n,\qquad
\text{with \
}\Psi_{(1)}=i^{n-1}F_n.
\end{equation}

  \item[\rm(c)]
Let $\mathbb F$ be a
real closed field
$\mathbb P$, let
$\mathbb P+\mathbb Pi$
$($see \eqref{1pp}$)$
be its algebraic
closure with
involution
$a+bi\mapsto a-bi$,
and let
$\varepsilon=\pm 1$.
One can take ${\cal
O}_{\mathbb F}$ to be
all $J_n(\lambda )$
with nonzero $\lambda
\in\mathbb P$ and all
realifications
$J_n(\lambda)^{\mathbb
P}$with
$\lambda\in(\mathbb
P+\mathbb
Pi)\smallsetminus\mathbb
P$ determined up to
replacement by
$\bar\lambda$.
\begin{itemize}
  \item[\rm(i)]
For $\lambda
\in\mathbb P$,
\begin{equation*}
J_n(\lambda
)_{(\varepsilon)}\
\text{ exists
}\quad\Longleftrightarrow\quad
\lambda =\pm 1\text{
and }
\varepsilon=(-1)^{n+1}.
\end{equation*}
If it exists then
$J_n(\lambda )$ is
similar to
\begin{equation*}\label{idq}
\Psi=\lambda
\Lambda_n,\qquad
\text{with \ }
\Psi_{(\varepsilon)}=F_n.
\end{equation*}
The field $\mathbb
T({\cal A}_{\Psi})$,
which is constructed
basing on the
endomorphisms of the
corresponding
selfadjoint
representation
\[\underline{\cal
A}_{\Psi}=(\lambda
\Lambda_n,F_n,
(\lambda
\Lambda_n)^*)\qquad
(\text{see
\eqref{nsp}}),\] is
naturally isomorphic
to $\mathbb P$.

  \item[\rm(ii)]
For
$\lambda,\mu\in(\mathbb
P+\mathbb
Pi)\smallsetminus\mathbb
P$,
\begin{equation*}\label{kust}
J_n(\lambda)^{\mathbb
P}\text{ is similar to
}(J_n(\mu)^{\mathbb
P})^{-T}
\quad\Longleftrightarrow\quad
\lambda\in\{{\mu}^{-1},
\bar{\mu}^{-1}\},
\end{equation*}
\begin{equation*}\label{ndsis}
J_n(\lambda)^{\mathbb
P}_{(\varepsilon)}\
\text{ exists
}\quad\Longleftrightarrow\quad
|\lambda |=1.
\end{equation*}
If it exists then
$J_n(\lambda)^{\mathbb
P}$ is similar to
\begin{equation}\label{ksihddy}
\Psi=(\lambda
\Lambda_n)^{\mathbb
P},\qquad \text{with \
}
 \Psi_{(\varepsilon)}
= (i^{n-(\varepsilon
+1)/2}F_n)^{\mathbb
P}.
\end{equation}
The field $\mathbb
T({\cal A}_{\Psi})$ is
naturally isomorphic
to $\mathbb P+\mathbb
Pi$ with involution
$a+bi\mapsto a-bi$.
\end{itemize}

  \item[\rm(d)]
Let $\mathbb F$ be the
skew field $\mathbb H$
of quaternions with
quaternionic
conjugation \eqref{ne}
or quaternionic
semiconjugation
\eqref{nen} over a
real closed field
$\mathbb P$, and let
$\varepsilon=\pm 1$.
One can take ${\cal
O}_{\mathbb F}$ to be
all $J_n(\lambda)$
with nonzero
$\lambda=a+bi\in\mathbb
P+\mathbb Pi$
determined up to
replacement by $a-bi$.
For nonzero
$\lambda,\mu\in\mathbb
P+\mathbb Pi$,
\begin{equation}\label{lust}
J_n(\lambda)\text{ is
similar to
}J_n(\mu)^{-*}
\quad\Longleftrightarrow\quad
\lambda\in\{{\mu}^{-1},
\bar{\mu}^{-1}\},
\end{equation}
\begin{equation*}\label{nsis}
J_n(\lambda)_{(\varepsilon)}\
\text{ exists
}\quad\Longleftrightarrow\quad
|\lambda |=1.
\end{equation*}

If it exists then
$J_n(\lambda)$ is
similar to
\begin{equation}\label{ksim1}
\Psi=\lambda\Lambda_n,\qquad
\text{with \ }
\Psi_{(\varepsilon
)}=i^{n-(\varepsilon
+1)/2}F_n,
\end{equation}
and
\begin{itemize}
  \item[\rm(i)]
if $\lambda\ne\pm 1$,
then the field
$\mathbb T({\cal
A}_{\Psi})$ is
naturally isomorphic
to $\mathbb P+\mathbb
Pi$ with involution
$a+bi\mapsto a-bi$,

  \item[\rm(ii)]
if $\lambda =\pm 1$
and
$\varepsilon=(-1)^{n+1}$,
then $\mathbb T({\cal
A}_{\Psi})$ is
naturally isomorphic
to $\mathbb F$ and
this isomorphism
preserves the
involution,

  \item[\rm(iii)]
if $\lambda =\pm 1$
and
$\varepsilon=(-1)^{n}$,
then $\mathbb T({\cal
A}_{\Psi})$ is
naturally isomorphic
to $\mathbb F$ and if
the involution on
$\mathbb F$ is
quaternionic
conjugation \eqref{ne}
or quaternionic
semiconjugation
\eqref{nen} then the
involution on $\mathbb
T({\cal A}_{\Psi})$ is
\eqref{nen} or
\eqref{ne},
respectively.
\end{itemize}
\end{itemize}
\end{lemma}

\begin{proof}
(a) Let $\mathbb F$ be
an algebraically
closed field with the
identity involution
and $\varepsilon=\pm
1$. By \eqref{lbdr}
and \eqref{4.adlw}, if
$J_n(\lambda)_{(\varepsilon
)}$ exists then
$\lambda=\pm 1$ and
$\varepsilon=(-1)^{n+1}$.
Let these conditions
be satisfied. Since
$\lambda\Lambda_n$ is
similar to
$\Psi=J_n(\lambda)$,
it remains to check
that
$\Psi_{(\varepsilon)}=F_n$
fulfils \eqref{mskhd},
that is,
\begin{equation}\label{mskhd1}
\Psi_{(\varepsilon)}
=\varepsilon
\Psi_{(\varepsilon)}^*,\qquad
\Psi_{(\varepsilon)}
=\Psi^*\Psi_{(\varepsilon)}\Psi.
\end{equation}

The first equality is
obvious. The second is
satisfied because
\begin{align*}
\Psi_{(\varepsilon)}^{-1}
&\cdot\Psi^*\Psi_{(\varepsilon)}
\Psi=F_n^{-1}\Lambda_n^T
F_n\Lambda_n\\
&=\begin{bmatrix}
0&&&1
\\&&-1&\\&1&&\\\ddd&&&0
\end{bmatrix}
\begin{bmatrix} 1&&&0
\\2&1&&\\
\vdots&\ddots&\ddots&\\
2&\dots&2&1
\end{bmatrix}
\begin{bmatrix}
0&&&\ddd
\\&&1&
\\&-1&\\
1&&&0
\end{bmatrix}\Lambda_n
\\&=\begin{bmatrix}
1&-2&2&\ddots&(-1)^{n-1}2
\\&1&-2&\ddots&\ddots
\\&&1&\ddots&2
\\&&&\ddots&-2
\\0&&&&1
\end{bmatrix}\begin{bmatrix}
1&2&\cdots&2
\\&1&\ddots&\vdots
\\&&\ddots&2
\\0&&&1
\end{bmatrix}\quad
(\text{write
}J:=J_n(0))
\\&=(I_n-2J+2J^2-2J^3+2J^4-\cdots)
(I_n+2J+2J^2+2J^3+\cdots)=I_n.
\end{align*}
\medskip

(b)  Let $\mathbb
F=\mathbb P+\mathbb
Pi$ be an
algebraically closed
field with nonidentity
involution represented
in the form
\eqref{1pp11}. By
\eqref{lbdr}, if
$J_n(\lambda)_{(1)}$
exists for $\lambda
=a+bi$ then $x-\lambda
=x-\bar\lambda^{-1}$.
Thus, $\lambda
=\bar\lambda^{-1}$ and
by \eqref{1ii}
$1=\lambda
\bar\lambda=a^2+b^2=|\lambda
|^2$.

Let $|\lambda|=1$. The
matrix
$\Psi=\lambda\Lambda_n$
in \eqref{ksim} is
similar to
$J_n(\lambda)$, the
first equality in
\eqref{mskhd1} is
obvious for
$\Psi_{(1)}=i^{n-1}F_n$
and the second holds
since it holds for
\eqref{ksidq}.
\medskip

(c) Let $\mathbb
F=\mathbb P$ be a real
closed field and
$\varepsilon=\pm 1$.
Let $\mathbb
K:=\mathbb P+\mathbb
Pi$ be the algebraic
closure of $\mathbb P$
represented in the
form \eqref{1pp}  with
involution
$a+bi\mapsto a-bi$. By
\cite[Theorem
3.4.5]{hor}, we can
take ${\cal
O}_{\mathbb P}$ to be
all $J_n(\lambda )$
with $0\ne\lambda
\in\mathbb P$ and all
$J_n(\lambda
)^{\mathbb P}$ with
$\lambda\in{\mathbb
K}\smallsetminus\mathbb
P$ determined up to
replacement by
$\bar\lambda$.

Let us consider
$J_n(\lambda )$ with
$\lambda \in\mathbb
P$. By \eqref{lbdr}
and \eqref{4.adlw}, if
$J_n(\lambda
)_{(\varepsilon )}$
exists then $\lambda
=\pm 1$ and
$\varepsilon=(-1)^{n+1}$.
Hence we can use
$\Psi$ and
$\Psi_{(\varepsilon)}$
from \eqref{ksidq}. In
view of \eqref{ksy}
and since
$p_{\Psi}(x)=x-\lambda
$,
$$\mathbb T({\cal
A}_{\Psi})\simeq
{\mathbb
P}(\kappa)={\mathbb
P}[x]/p_{\Psi}(x){\mathbb
P}[x]\simeq {\mathbb
P}.$$

Let now
$\Phi:=J_n(\lambda
)^{\mathbb P}$ with
$\lambda \in{\mathbb
K}\smallsetminus\mathbb
P$.  Then
\begin{equation}\label{kgv}
p_{\Phi}(x)=(x-\lambda)(x-\bar\lambda)
=x^2-(\lambda+\bar\lambda)
+|\lambda|^2.
\end{equation}
If $J_n(\lambda
)^{\mathbb
P}_{(\varepsilon )}$
exists, then
$|\lambda|=1$ by
\eqref{lbdr} and
\eqref{kgv}.

If $\lambda
\in{\mathbb
K}\smallsetminus\mathbb
P$ and $|\lambda|=1$,
then we can take
$\Psi$ and
$\Psi_{(\varepsilon)}$
as in \eqref{ksihddy}
due to the following
observation. The
equalities
\eqref{mskhd1} hold
for \eqref{ksim1}
since for $\varepsilon
=1$ they were checked
in (b) and so they are
fulfilled for
$\varepsilon =-1$ too.
Therefore,
\eqref{mskhd1} hold
true for
\eqref{ksihddy} due to
the following property
of realification:
\emph{if matrices
$M_1,\dots,M_l,M_1^*,\dots,M_l^*$
over $\mathbb
P+\mathbb Pi$ satisfy
an equation with
coefficients in
$\mathbb P$, then
their realifications
also satisfy the same
equation}. This
property is valid
since for each matrix
$M=A+Bi$ with $A$ and
$B$ over $\mathbb P$,
its realification
$M^{\mathbb P}$ (see
\eqref{1j}), up to
simultaneous
permutations of rows
and columns, has the
form
\begin{equation}\label{luf}
M_{\mathbb
P}:=\begin{bmatrix}
  A&-B\\B&A
\end{bmatrix}=S^{-1}(M\oplus \bar
M)S=S^{*}(M\oplus \bar
M)S
\end{equation}
with
\[
S:=\frac{1}{\sqrt{2}}
\begin{bmatrix}
  I&iI\\I&-iI
\end{bmatrix}=S^{-*};
\]
the middle equality in
\eqref{luf} follows
from
\[
\begin{bmatrix}
  A+Bi&0\\0&A-Bi
\end{bmatrix}
\begin{bmatrix}
  I&iI\\I&-iI
\end{bmatrix}=
\begin{bmatrix}
  I&iI\\I&-iI
\end{bmatrix}
\begin{bmatrix}
  A&-B\\B&A
\end{bmatrix}.
\]

By \eqref{ksy} and
since $\deg
p_{\Psi}(x)=2$,
$\mathbb T({\cal
A}_{\Psi})\simeq{\mathbb
P}(\kappa)\simeq\mathbb
P+\mathbb Pi.$ The
involution \eqref{kxu}
on $\mathbb T({\cal
A}_{\Psi})$ is not the
identity (otherwise
$\kappa= \kappa^{-1}$,
$\kappa^2-1=0$, i.e.
$p_{\Psi}(x)=x^2-1$,
which contradicts the
irreducibility of
$p_{\Psi}(x)$).
\medskip

(d) Let $\mathbb F$ be
the skew field
$\mathbb H$ of
quaternions over a
real closed field
$\mathbb P$, and let
$\varepsilon=\pm 1$.
By \cite[Section 3, \S
12]{jac}, we can take
${\cal O}_{\mathbb F}$
to be all
$J_n(\lambda)$ which
$\lambda=a+bi\in\mathbb
P+\mathbb Pi$
determined up to
replacement by $a-bi$.
For any nonzero
$\mu\in\mathbb
P+\mathbb Pi$, the
matrix $J_n(\mu)^{-*}$
is similar to
$J_n(\bar\mu^{-1})$,
by (b) it is similar
to $J_n(\lambda)$ with
$\lambda \in\mathbb
P+\mathbb Pi$ if and
only if
$\lambda\in\{{\mu}^{-1},
\bar{\mu}^{-1}\}$.
This proves
\eqref{lust}.

 Using
\eqref{lbdr} and
reasoning as in (b),
we make sure that if
$J_n(a+bi)_{(\varepsilon)}$
exists then $a^2+b^2=
1$. We can take $\Psi$
and
$\Psi_{(\varepsilon
)}$ as in
\eqref{ksim1} since
the equalities
\eqref{mskhd1} for
them were checked in
(c).

Due to \eqref{mdtc},
$[G_1,G_2]\in \End
(\underline{\cal
A}_{\Psi})$ if and
only if
\begin{equation}\label{heu}
G_1\Psi=\Psi
G_1,\qquad
G_2\Psi_{(\varepsilon)}=
\Psi_{(\varepsilon)}
G_1,\qquad
G_2\Psi^*=\Psi^* G_2.
\end{equation}
The last equality
follows from the
others:
\begin{align*}
\Psi^* G_2 &= \Psi^*
\Psi_{(\varepsilon)}
G_1
\Psi_{(\varepsilon)}^{-1}
= \Psi_{(\varepsilon)}
\Psi^{-1}
G_1 \Psi_{(\varepsilon)}^{-1}\\
& =
\Psi_{(\varepsilon)}G_1
\Psi^{-1}
\Psi_{(\varepsilon)}^{-1}=
\Psi_{(\varepsilon)}G_1
\Psi_{(\varepsilon)}^{-1}
\Psi^*= G_2\Psi^*.
\end{align*}

Let $\lambda\in\mathbb
K:=\mathbb P+\mathbb
Pi$ and $|\lambda|=1$.

(i) First we consider
the case $\lambda\ne
\pm 1$. Represent
$G_1$ in the form
$U+Vj$, where
$U,V\in\mathbb
K^{n\times n}$. Then
the first equality in
\eqref{heu} becomes
$
(U+Vj)\lambda\Lambda
=\lambda\Lambda(U+Vj)
$
and falls into two
equalities
\[
U\lambda\Lambda
=\lambda\Lambda
U,\qquad
V\bar\lambda\Lambda j
=\lambda\Lambda Vj
\]
(quaternionic
conjugation \eqref{ne}
and quaternionic
semiconjugation
\eqref{nen} coincide
on $\mathbb K$). By
the second equality,
\[
(\bar\lambda-\lambda)V=
\lambda
(\Lambda-I)V-\bar\lambda
V(\Lambda-I)
\]
and so $V=0$ since
$\lambda\ne\bar\lambda$
and because
$\Lambda-I$ is
nilpotent upper
triangular. By the
first equality (which
is over the field
$\mathbb K$),
$G_1=U=f(\lambda
\Lambda)=f(\Psi)$ for
some $f\in\mathbb
K[x]$; see the
beginning of the proof
of Lemma
\ref{lenhi}(d). Since
$\Psi_{(\varepsilon)}$
is over $\mathbb K$
and in view of
\eqref{heu},
\[G_2=\Psi_{(\varepsilon)}G_1
\Psi_{(\varepsilon)}^{-1}
=
f(\Psi_{(\varepsilon)}\Psi
\Psi_{(\varepsilon)}^{-1})
=f(\Psi^{-*}).\]

Due to
$G_2=\Psi_{(\varepsilon)}G_1
\Psi_{(\varepsilon)}^{-1}$,
the homomorphism
$[G_1,G_2]\in \End
(\underline{\cal
A}_{\Psi})$ is
completely determined
by $G_1=f(\Psi)$. The
matrix $\Psi=\lambda
\Lambda$ is upper
triangular, so the
mapping
$f(\Psi)\mapsto
f(\lambda)$,
$f\in\mathbb K[x]$,
defines an
endomorphism of rings
$\End (\underline{\cal
A}_{\Psi})\to \mathbb
K$, its kernel is the
radical of $\End
(\underline{\cal
A}_{\Psi})$. Hence
$\mathbb T( {\cal
A}_{\Psi})$ can be
identified with
$\mathbb K$. In view
of \eqref{ami}, the
involution on $\mathbb
T({\cal A}_{\Psi})$ is
induced by the mapping
$G_1\mapsto G_2^*$ of
the form
\[f(\lambda
\Lambda)\ \longmapsto\
f((\lambda
\Lambda)^{-*})^*= \bar
f((\lambda
\Lambda)^{-1}).\]
Therefore, the
involution is
\[f(\lambda)\ \longmapsto\  \bar
f({\lambda}^{-1})=
\bar f(\bar{\lambda})=
\overline{f({\lambda})}\]
and coincides with the
involution
$a+bi\mapsto a-bi$.
\medskip

(ii)\&(iii) Let
$\lambda= \pm 1$.
Define
\begin{align*}
\check
h&:=a+bi-cj-dk\quad
\text{for each}\ \
h=a+bi+cj+dk\in\mathbb
H,
\\
\check f(x)&:=\sum_l
\check h_lx^l\quad
\text{for each}\ \
f(x)=\sum_l
h_{l}x^l\in \mathbb
H[x].
\end{align*}

Because $\lambda= \pm
1$ and by the first
equality in
\eqref{heu}, $G_1$ has
the form
\[
G_1=\begin{bmatrix}
 a_1&a_2&\ddots&a_{n}
 \\&a_1&\ddots&\ddots
 \\&&\ddots&a_2\\
 0&&&a_1
\end{bmatrix}, \qquad
a_1,\dots,a_n\in\mathbb
H.
\]
Thus, $G_1=f(\Psi)$
for some
$f(x)\in\mathbb H[x]$.

Using the second
equalities in
\eqref{heu} and
\eqref{ksim1} and the
identity $if(x)=\check
f(ix)$, we obtain
\begin{align*}
G_2&=\Psi_{(\varepsilon)}
G_1
\Psi_{(\varepsilon)}^{-1}
=\Psi_{(\varepsilon)}
f(\Psi)
\Psi_{(\varepsilon)}^{-1}\\
&=\begin{cases}
   f( \Psi_{(\varepsilon)}
\Psi
\Psi_{(\varepsilon)}^{-1})=
f( \Psi^{-*}) &
\text{if $
    \varepsilon=(-1)^{n+1}$,} \\
    \check f( \Psi_{(\varepsilon)}
\Psi
\Psi_{(\varepsilon)}^{-1})=
\check f( \Psi^{-*}) &
\text{if $
    \varepsilon=(-1)^{n}$}.
  \end{cases}
\end{align*}

Since the homomorphism
$[G_1,G_2]$ is
completely determined
by $G_1=f(\Psi)$, the
matrix $\Psi=\lambda
\Lambda_n$ is upper
triangular, its main
diagonal is
$(\lambda,\dots,\lambda)$,
and $\lambda=\pm 1$,
we conclude that the
mapping
$f(\Psi)\mapsto
f(\lambda)$ defines an
endomorphism of rings
$\End(\underline{\cal
A}_{\Psi})\to \mathbb
H$ and its kernel is
the radical of
$\End(\underline{\cal
A}_{\Psi})$. Hence
$\mathbb T( {\cal
A}_{\Psi})$ can be
identified with
$\mathbb H$. The
involution on $\mathbb
T( {\cal A}_{\Psi})$
is induced by the
mapping $G_1\mapsto
G_2^*$, that is, by
\[
f(\Psi)\mapsto
\begin{cases}
   \bar{f}( \Psi^{-1}) &
\text{if $
    \varepsilon=(-1)^{n+1}$,} \\
   \widehat{f}(
\Psi^{-1}) & \text{if
$
\varepsilon=(-1)^{n}$}.
  \end{cases}
\]
Here $h\mapsto\bar h$
is the involution on
$\mathbb F$ that is
quaternionic
conjugation \eqref{ne}
or quaternionic
semiconjugation
\eqref{nen}, and
$h\mapsto\widehat{h}$
denotes the remaining
involution \eqref{nen}
or \eqref{ne}. Thus,
the involution on
$\mathbb T( {\cal
A}_{\Psi})$ is
$h\mapsto\bar h$ if
$\varepsilon=(-1)^{n+1}$
and $h\mapsto
\widehat{h}$ if
$\varepsilon=(-1)^{n}$.
\end{proof}

\begin{proof}[Proof of
Theorem \ref{theor}]
Let $\mathbb F$ be one
of the fields and skew
fields considered in
Theorem \ref{theor}.
By Theorem
\ref{tetete}, each
representation of a
dograph $D$ over
$\mathbb F$ is
uniquely, up to
isomorphism of
summands, decomposes
into a direct sum of
indecomposable
representations. Hence
the problem of
classifying its
representations
reduces to the problem
of classifying
indecomposable
representations.

Let ${\cal O}_{\mathbb
F}$ be a maximal set
of nonsingular
indecomposable
canonical matrices
over $\mathbb F$ for
similarity. Due to
Theorem \ref{tetete}
and Lemma \ref{lenhi},
the following
representations form a
maximal set of
nonisomorphic
indecomposable
representations of the
dograph $D$ defined in
\eqref{jsos+}:
\begin{itemize}
  \item[(i)]
${\cal
M}_{\Phi}^+=(\Phi\oplus
\Phi^{-*}, I\dia
\varepsilon I)$, in
which $\Phi\in{\cal
O}_{\mathbb F}$ is
such that
$\Phi_{(\varepsilon)}$
does not exist; $\Phi$
is determined up to
replacement by the
matrix $\Psi\in{\cal
O}_{\mathbb F}$ that
is similar to
$\Phi^{-*}$.

  \item[(ii)]
${\cal
A}_{\Phi}=(\Phi,
\Phi_{(\varepsilon)})$
and ${\cal
A}_{\Phi}^-=(\Phi,
-\Phi_{(\varepsilon)})$,
in which $\Phi\in{\cal
O}_{\mathbb F}$ is
such that
$\Phi_{(\varepsilon)}$
exists. The
representation ${\cal
A}_{\Phi}^-$ is
withdrawn if ${\cal
A}_{\Phi}\simeq {\cal
A}_{\Phi}^-$, this
occurs if and only if
\begin{itemize}
  \item
$\mathbb F$ is an
algebraically closed
field with the
identity involution,
or
  \item
$\mathbb F$ is not an
algebraically closed
field and either
$\mathbb T({\cal
A}_{\Phi})$ is an
algebraically closed
field with the
identity involution,
or $\mathbb T({\cal
A}_{\Phi})$ is a skew
field of quaternions
with involution
different from
quaternionic
conjugation
\eqref{ne}.
\end{itemize}
\end{itemize}

Thus, the statements
(a)--(d) of Theorem
\ref{theor} follow
from the statements
(a)--(d) of Lemma
\ref{lsdy}.
\end{proof}

\section{Metric and
selfadjoint operators
with respect to
degenerate forms}
\label{metric}

Recall that a
classification problem
is called \emph{wild}
if it contains the
problem of classifying
pairs of matrices up
to simultaneous
similarity and hence
it contains the
problem of classifying
representations of
each quiver. A linear
operator is called
\emph{metric} or
\emph{selfadjoint}
with respect to a form
$B$ (possibly,
degenerate) if $
B(Au,Av)=B(u,v)$ or
$B(Au,v)=B(u,Av)$,
respectively, for all
$u$ and $v$. The
following theorem was
proved in
\cite[Theorem
5.4]{ser_prep}.

\begin{theorem}\label{th_dic}
The problem of
classifying pairs
$(A,B)$ consisting of
a form $B$ on a vector
space $V$ over a field
of characteristic
different from $2$ and
an operator $A$ that
is metric with respect
to $B$ is wild in each
of the following three
cases: $B$ is
symmetric, $B$ is
skew-symmetric, or $B$
is Hermitian. This
statement also holds
if `metric'' is
replaced by
``selfadjoint''.
\end{theorem}

\begin{proof}
(a) Suppose first that
$A$ is metric. The
problem of classifying
pairs $(A,B)$ is given
by the dograph
\eqref{jsos+} without
the condition
``$\beta$ is
nonsingular'' and
reduces to the problem
of classifying
representations of the
corresponding quiver
\eqref{jso+} without
the condition
``$\beta$ is
bijective''. Each
representation of this
quiver has the form
\eqref{ser14sd}
without the condition
``$B$ is
nonsingular'', i.e.,
it is given by
matrices $A$, $B$, $C$
of sizes $m\times m$,
$n\times m$, $n\times
n$ satisfying the
relation:
\begin{equation}\label{sev}
B=CBA.
\end{equation}
By \eqref{msp}, two
matrix triples
$(A,B,C)$ and
$(A',B',C')$ give
isomorphic
representations if and
only if there exist
nonsingular matrices
$R$ and $S$ such that
\begin{equation}\label{pfs}
RA=A'R,\quad
SB=B'R,\quad SC=C'S.
\end{equation}
Since
$B'=SBR^{-1}=I_r\oplus
0$ for some
nonsingular $R$ and
$S$, it suffices to
consider only the
triples $(A,B,C)$ with
$B=I_r\oplus 0$. Such
a triple satisfies
\eqref{sev} if and
only if it has the
form
\begin{equation}\label{grx}
(A,B,C)=\left(\begin{bmatrix}
A_{11}&0\\
 A_{21}&A_{22}
\end{bmatrix},\ \begin{bmatrix}
I_r&0\\
 0&0
\end{bmatrix},\
\begin{bmatrix}
C_{11}&C_{12}\\
0&C_{22}
\end{bmatrix}\right),
\end{equation}
in which $A_{11}$ and
$C_{11}$ are $r\times
r$ matrices and
$C_{11}A_{11}=I_r$.
Triples \eqref{grx}
and
\begin{equation}\label{grxs}
(A',B',C')=\left(\begin{bmatrix}
A_{11}'&0\\
 A_{21}'&A_{22}'
\end{bmatrix},\ \begin{bmatrix}
I_r&0\\
 0&0
\end{bmatrix},\
\begin{bmatrix}
C_{11}'&C_{12}'\\
0&C_{22}'
\end{bmatrix}\right)
\end{equation}
give isomorphic
representations if and
only if there exist
nonsingular $R$ and
$S$ satisfying
\eqref{pfs}. The
equality $SB=B'R$ with
$B'=B=I_r\oplus 0$
implies
\begin{equation}\label{suc}
R=\begin{bmatrix}
R_{11}&0\\
R_{21}&R_{22}
\end{bmatrix},\qquad
S=\begin{bmatrix}
S_{11}&S_{12}\\
0&S_{22}
\end{bmatrix},\qquad
R_{11}=S_{11}.
\end{equation}
The remaining
equalities in
\eqref{pfs} take the
form
\begin{equation}\label{msyj}
\begin{bmatrix}
R_{11}&0\\
R_{21}&R_{22}
\end{bmatrix}
\begin{bmatrix}
A_{11}&0\\
 A_{21}&A_{22}
\end{bmatrix}
=
\begin{bmatrix}
A_{11}'&0\\
 A_{21}'&A_{22}'
\end{bmatrix}
\begin{bmatrix}
R_{11}&0\\
R_{21}&R_{22}
\end{bmatrix},
\end{equation}
\begin{equation}\label{msy}
\begin{bmatrix}
S_{11}&S_{12}\\
0&S_{22}
\end{bmatrix}
\begin{bmatrix}
C_{11}&C_{12}\\
0&C_{22}
\end{bmatrix}
=
\begin{bmatrix}
C_{11}'&C_{12}'\\
0&C_{22}'
\end{bmatrix}
\begin{bmatrix}
S_{11}&S_{12}\\
0&S_{22}
\end{bmatrix}.
\end{equation}

Therefore, the problem
of classifying pairs
$(A,B)$ contains the
problem of classifying
upper block-triangular
matrices
\begin{equation}\label{lkd}
\begin{bmatrix}
C_{11}&C_{12}\\
0&C_{22}
\end{bmatrix},
\end{equation}
in which $C_{11}$ is
nonsingular, with
respect to upper
block-triangular
similarity. The
wildness of this
problem and many
analogous problems was
proved, for example,
in
\cite[Section
3.3.1]{ser-can}.
\medskip

(b) Suppose now that
$A$ is selfadjoint.
The problem of
classifying pairs
$(A,B)$ is given by
the dograph
\eqref{jsos+} in
which all the
relations are replaced
by
\begin{equation}
\label{mdr}
 \beta=\varepsilon \beta^*,\quad
  \beta\alpha=\alpha^*\beta.
\end{equation}
It reduces to the
problem of classifying
representations of the
corresponding quiver
\eqref{jso+} with
relations \eqref{mdr}.
Each representation of
this quiver is given
by matrices $A$, $B$,
$C$ of sizes $m\times
m$, $n\times m$,
$n\times n$ such that
\begin{equation}\label{ssev}
BA=CB.
\end{equation}
Let us consider
triples $(A,B,C)$ with
$B=I_r\oplus 0$. Such
a triple satisfies
\eqref{ssev} if and
only if it has the
form \eqref{grx} with
$C_{11}=A_{11}$.
Triples \eqref{grx}
and \eqref{grxs} give
isomorphic
representations if and
only if the equalities
\eqref{msyj} and
\eqref{msy} are valid
for some nonsingular
$R$ and $S$ of the
form \eqref{suc}.

Therefore, the problem
of classifying pairs
$(A,B)$ contains the
wild problem of
classifying upper
block-triangular
matrices \eqref{lkd}
with respect to upper
block-triangular
similarity.
\end{proof}


\begin{thebibliography}{99}
\parskip=-1.5pt


\bibitem{apl-Rod}
D. Alpay, A.C.M. Ran,
L. Rodman, Basic
classes of matrices
with respect to
quaternionic
indefinite inner
product spaces,
\emph{Linear Algebra
Appl.} 416 (2006)
242--269.

\bibitem{au-Rod}
Y.-H. Au-Yeung, C.-K.
Li, L. Rodman,
H-unitary and Lorentz
matrices: A review,
\emph{SIAM J. of
Matrix Analysis}, 25
(2004) 1140--1162.

\bibitem{bel-ser_complexity}
G.R. Belitskii, V.V.
Sergeichuk, Complexity
of matrix problems,
\emph{Linear Algebra
Appl.} 361 (2003)
203--222.

\bibitem{bour2}
N. Bourbaki,
\emph{Alg\`{e}bre,
Chaps. 8,9},
Actualiti\'{e}s Sci.
Indust., nos. 1261,
1272, Hermann, Paris,
1958, 1959.

\bibitem{der-wey}
H. Derksen, J. Weyman,
Generalized quivers
associated to
reductive groups,
\emph{Colloq. Math.}
94 (2002) 151--173.

\bibitem{gab}
P. Gabriel,
Unzerlegbare
Darstellungen I, {\it
Manuscripta Math.} 6
(1972) 71--103.

\bibitem{gab2}
P. Gabriel, Appendix:
degenerate bilinear
forms, {\it J.
Algebra} 31 (1974)
67--72.

\bibitem{gab_vos}
P. Gabriel, L.A.
Nazarova, A.V. Roiter,
V.V. Sergeichuk, D.
Vossieck, Tame and
wild subspace
problems, {\it
Ukrainian Math. J.} 45
(1993) 335--372.

\bibitem{g-l-r}
I. Gohberg, P.
Lancaster, L. Rodman,
\emph{Matrices and
indefinite scalar
products}, Operator
theory: advances and
applications, vol. 8,
BirkhÄauser Verlag,
Basel, 1983.

\bibitem{g-l-r1}
I. Gohberg, P.
Lancaster, L. Rodman,
\emph{Indefinite
linear algebra and
applications},
Birkh\"auser Verlag,
Basel, 2005.

\bibitem{goh}
I. Gohberg, B.
Reichstein, On
classification of
normal matrices in an
indefinite scalar
product, {\it Integral
Equations and Operator
Theory} 13 (1990)
364--394.

\bibitem{hor}
R.A. Horn, C.R.
Johnson, {\it Matrix
Analysis}, Cambridge
U. P., Cambridge,
1985.

\bibitem{hor-ser}
R.A. Horn, V.V.
Sergeichuk, Congruence
of a square matrix and
its transpose,
\emph{Linear Algebra
Appl.} 389 (2004)
347--353.

\bibitem{hor-ser_regul}
R.A. Horn, V.V.
Sergeichuk, A
regularization
algorithm for matrices
of bilinear and
sesquilinear forms,
\emph{Linear Algebra
Appl.} 412 (2006)
380--395.

\bibitem{hor-ser_can}
R.A. Horn, V.V.
Sergeichuk, Canonical
forms for complex
matrix congruence and
*congruence,
\emph{Linear Algebra
Appl.} 416 (2006)
1010--1032.

\bibitem{huppert1}
B. Huppert, Isometrien
von Vektorr\"{a}umen.
I, {\it Arch. Math.}
(Basel) 35 (1980)
164--176.

\bibitem{huppert2}
B. Huppert, Isometrien
von Vektorr\"{a}umen.
II, {\it Math. Z.} 175
(1980) 5--20.


\bibitem{jac}
N. Jacobson, {\it The
theory of rings},
Mathematical Surveys,
2. Providence, R. I.:
American Mathematical
Society, 1966.

\bibitem{kac}
V.G. Kac, Infinite
root systems,
representation of
graphs and invariant
theory, \emph{Invent.
Math.} 56 (1980)
57--92.

\bibitem{kar}
M. Karow, Self-adjoint
operators and pairs of
Hermitian forms over
the quaternions,
\emph{Linear Algebra
Appl.} 299 (1999)
101--117.

\bibitem{len}
S. Lang, {\it
Algebra}, Reading, MA:
Addison-Wesley, 1965.

\bibitem{mehl-rod_deg}
C. Mehl, L. Rodman,
Symmetric matrices
with respect to
sesquilinear forms,
\emph{Linear Algebra
Appl.} 349 (2002)
55--75.

\bibitem{miln}
J. Milnor, On
isometries of inner
product spaces, {\it
Invent. Math.} 8
(1969) 83--97.

\bibitem{naz_bon}
L.A. Nazarova, A.V.
Roiter, V.V.
Sergeichuk, V.M.
Bondarenko,
Application of modules
over a dyad for the
classification of
finite $p$-groups
possessing an abelian
subgroup of index $p$
and  of pairs of
mutually annihilating
operators,  \emph{Zap.
Nau\v cn. Sem.
Leningrad. Otdel. Mat.
Inst. Steklov. (LOMI)}
28 (1972), 69--92;
translation in {\it J.
Soviet Math.} 3 (no.
5) (1975) 636--654.

\bibitem{que-scha}
H.-G. Quebbermann, W.
Scharlau, M. Schulte,
Quadratic and
Hermitian forms in
additive and abelian
categories, {\it J.
Algebra} 59 (1979)
264--289.

\bibitem{riehm1}
C. Riehm, The
equivalence of
bilinear forms, {\it
J. Algebra} 31 (1974)
45--66.

\bibitem{riehm2}
C. Riehm, M.
Shrader-Frechette, The
equivalence of
sesquilinear forms,
{\it J. Algebra} 42
(1976) 495--530.


\bibitem{rod_pert}
L. Rodman, Similarity
vs unitary similarity
and perturbation
analysis of sign
characteristic:
Complex and real
indefinite inner
products, \emph{Linear
Algebra Appl.} 416
(2006) 945--1009.

\bibitem{rod_pair_nonst}
L. Rodman, Canonical
forms for symmetric
quaternionic matrix
pencils, submitted for
publication.

\bibitem{rod_pair_stand}
L. Rodman, Canonical
forms for mixed
symmetric/skew-symmetric
quaternionic matrix
pencils, submitted
for publication.

\bibitem{rod_pair_herm}
L. Rodman, Pairs of
hermitian and
skew-hermitian
quaternionic matrices:
canonical forms and
their applications, in
preparation.

\bibitem{roi}
A.V. Roiter, Bocses
with involution, in:
{\it Representations
and quadratic forms},
Akad. Nauk Ukrain.
SSR, Inst. Mat., Kiev,
1979, 124--128 (in
Russian).

\bibitem{sch}
R. Scharlau, Zur
Klassification von
Bilineaformen und von
Isometrien \"{u}ber
K\"{o}rpern, {\it
Math. Z.} 178 (1981)
359--373.

\bibitem{w.schar}
W. Scharlau,
{\it Quadratic and
Hermitian Forms},
Springer-Verlag, 1985.

\bibitem{ser_first}
V.V. Sergeichuk,
Representations of
simple involutive
quivers, in:
\emph{Representations
and quadratic forms},
Akad. Nauk Ukrain.
SSR, Inst. Mat., Kiev,
1979, 127--148 (in
Russian).

\bibitem{ser_disch}
V.V. Sergeichuk,
Representation of
oriented schemes, in:
\emph{Linear algebra
and the theory of
representations},
Akad. Nauk Ukrain.
SSR, Inst. Mat., Kiev,
1983, 110--134 (in
Russian).

\bibitem{ser_prep}
V.V. Sergeichuk, {\it
Classification
problems for systems
of linear mappings and
sesquilinear forms},
Preprint, Kiev
University, 1983, 60
p. (in Russian) =
Manuscript No. 196
Uk-D84, deposited at
the Ukrainian NIINTI,
1984; R. Zh. Mat.
1984, 7A331.

\bibitem{ser_izv}
V.V. Sergeichuk,
Classification
problems  for systems
of forms and linear
mappings, {\it  Math.
USSR-Izv.} 31 (no. 3)
(1988) 481--501.

\bibitem{ser_sub}
V.V. Sergeichuk,
Classification of
pairs of subspaces in
scalar product spaces,
{\it Ukrainian Math.
J.} 42 (no. 4) (1990)
487--491.

\bibitem{ser1}
V.V. Sergeichuk,
Classification of
sesquilinear forms,
pairs of hermitian
forms, self-conjugate
and isometric
operators over the
division ring of
quaternions, {\it
Math. Notes} 49 (no.
3--4) (1991) 409--414.

\bibitem{ser_sym}
V.V. Sergeichuk,
Symmetric
representations of
algebras with
involution,
\emph{Math. Notes} 50
(no. 3-4) (1991)
1058--1061.

\bibitem{ser_unit}
V.V. Sergeichuk,
Unitary and Euclidean
representations of a
quiver, \emph{Linear
Algebra Appl.} 278
(1998) 37--62.

\bibitem{ser-can}
V.V. Sergeichuk,
Canonical matrices for
linear matrix
problems, {\it Linear
Algebra Appl.} 317
(2000) 53--102.

\bibitem{shme}
D.A. Shmelkin, Signed
quivers, symmetric
quivers and root
systems, \emph{J.
London Math. Soc.} (2)
73 (2006) 586--606.

\bibitem{wan}
B.L. van der Waerden,
{\it Algebra}. New
York, NY:
Springer-Verlag, 1991.
\end{thebibliography}
\end{document}